\documentclass[reqno, 10.5pt]{amsart}
\usepackage{amsfonts, amssymb, amsthm, amsmath}
\usepackage{url}
\usepackage{epsfig}
\usepackage{tikz}
\usetikzlibrary{positioning}
\numberwithin{equation}{section}
\theoremstyle{definition}

\newtheorem{cor}{Corollary}[section]
\newtheorem{defn}{Definition}[section]
\newtheorem{examp}{Example}[section]
\newtheorem{lem}{Lemma}[section]
\newtheorem{prop}{Proposition}[section]
\newtheorem{rem}{Remark}[section]

\begin{document}

% Title of the paper
%\title{The KNW Quick Simulation Random Fields}
\title{Properties of Quick Simulation Random Fields}
% Authors of the paper

\author{Biao Wu}
\email{biao.wu@ualberta.ca}

\author{Michael A. Kouritzin*}
\address{Department of Mathematical and Statistical Sciences\\
University of Alberta\\
Edmonton, Alberta\\
T6G 2G1 CANADA}
\email{mkouritz@math.ualberta.ca}
\thanks{*The author gratefully acknowledges support from NSERC through a Discovery Grant.}

\author{Fraser Newton}
\email{fnewton@math.ualberta.ca}

% Keywords

\maketitle

% Abstract of your paper
\begin{abstract}
Herein, we introduce and study a new class of discrete random fields designed for quick simulation and covariance inference under inhomogeneous condition.
Simulation of these correlated fields can be done in a single pass instead of relying on multi-pass convergent methods like the Gibbs Sampler or other Markov Chain
Monte Carlo methods. The fields are constructed directly from specified marginal probability mass functions and covariances between nearby sites.
The proposition on which the construction is based establishes when and how it is possible to simplify the conditional probabilities of each site
given the other sites in a manner that makes simulation quite feasible yet maintains desired marginal probabilities and covariances between sites.
Special cases of these correlated fields have been deployed successfully in data authentication, object detection and image generation. The limitations that must be imposed on the covariances and marginal probabilities in order for the algorithm to work are studied. What's more, a necessary and sufficient condition that guarantees the permutation property of correlated random fields are investigated. In particular, Markov random fields as a subclass of correlated random fields  are derived by a general and natural condition. Consequently, a direct and flexible single pass algorithm for simulating Markov random fields follows.
\end{abstract}

\section{Introduction}
\label{sec:Introduction}
\setcounter{subsection}{1}

Random fields are widely used in sciences and technologies to model spatially distributed random phenomena or objects. Within sciences, random fields are used in geophysics, astrophysics, statistical mechanics, underwater acoustics, structural biology and agriculture. Applications of random fields in technologies include TV signal processing, image processing in photography such as medical images (human brain imaging, functional magnetic resonance imaging, mammography), computer vision, web data extraction, clustering gene expression time series, natural language processing etc. Readers are referred to \cite{Ashburner03}, \cite{Chellappa93}, \cite{LiChang08}, \cite{LiHD95}, \cite{Li}, \cite{Winkler}, \cite{Worsley95}, \cite{Zhang01}, and \cite{Zhu}  for those applications. Technologically, researchers of random fields have dealt either with the modeling of images (for synthesis, recognition or compression purposes) or with the resolution of various spatial inverse problems (image restoration and reconstruction, deblurring, classification, segmentation, data fusion, optical flow estimation, optical character recognition, stereo matching, finger print classification, pattern recognition, face recognition, intelligent video surveillance, sparse signal recovery, natural language processing like Chinese chunk and so on, see \cite{Blue93evaluationof}, \cite{Chellappa95}, \cite{Li}, \cite{Sun2008}, and \cite{Winkler}).

Scientists and technicians are interested in the inverse problems such as image restoration, boundary detection,
tomographic reconstruction, shape detection from shading, and motion analysis. Many precisely formulated mathematical
models were constructed to model certain types of random fields, and various methods and estimators have been developed
to make the proposed models work in application. There are diverse needs calling for simulating random fields.
For example, simulation is employed to calculate minimum mean square (MMS) and maximum
posterior marginal (MPM) estimators, see \cite{Winkler}.
Simulation can also be a potential smoothing technique.
In the chapter 2 of Winkler \cite{Winkler}, various smoothing techniques were proposed
to clean ``dirty" pictures. Most of these methods involve simulation.
The difficult problem is how do we effectively simulate random fields.
A typical simulation would involve 100,000 or more highly correlated random variables,
which would certainly exceed the capacity of modern computers if one tried to
simulate the whole random field directly.

Researchers frequently resort to imposing discrete Markov assumptions on their random fields
to be simulated out of practical need.
In this regard, the Gibbs sampler was proposed to ease this simulation difficulty.
Briefly speaking, a Gibbs sampler starts with a given \emph{initial configuration}
(i.e. potential realization of the random field) or a configuration chosen at random
from some initial distribution, and then updates its configuration site by site based on the local characteristics of the random field. Once all sites of a configuration are sequentially updated, a \textbf{sweep} or a \textbf{pass} is finished. A Gibbs sampler usually takes hundreds of sweeps to produce a configuration closely consistent with a given distribution and there are still computational and convergence issues to deal with.

In this paper, we propose a new class of discrete \textbf{correlated random fields} which incorporate given probability mass functions (pmfs) \(\{\pi_{s_i}\}\) for all sites \(S=\{s_i\}_{i=1}^N\) and given covariances \(\beta_{s_i,s_j}\) between nearby sites. These fields are dsigned with efficient simulation in mind. The number of possible random configurations within a general discrete random field can be enormous and simulation is further complicated when the sites are correlated with one another.
These factors can make Gibbs sampling and other Markov chain Monte Carlo simulation impractical. However, Proposition \ref{mainprop:alg} on which our fields are based establishes a method to imbed desired covariances and marginal probabilities into
a random field while maintaining simulation ease. Indeed, Proposition \ref{mainprop:alg} is a simple means to construct \emph{some}
site-by-site conditional probabilities consistent with given marginal
probabilities and site-to-site covariances in such a way that sampling
the missing portion of a random field sequentially is very feasible.
More precisely, when simulating a new site, we compute this conditional probability mass
function of its state conditioned on the known portion and the previously-simulated sites.
This construction establishes the \emph{sequential simulation property} of our correlated random field,
that is to say, we can actually construct a random field in one pass based on this algorithm,
reducing the computation over the Gibbs sampler dramatically. For demonstration purposes, we discuss application of our random fields and simulation
algorithm to Data Authentication, Object Detection and Image Creation.

The constraints and properties of the random fields generated by Proposition \ref{mainprop:alg} are discussed in detail. The necessary and sufficient conditions of the regularity type are given in Proposition \ref{prop:condition} when base set is a singleton. In particular, we investigate the conditions related to the marginality and permutation properties of our random fields. In the case where one wants to match the covariances between each pair of sites for \(S=\{s_i\}_{i=1}^n\), the field takes the form (See Lemma \ref{lem:jointpmf})
\begin{eqnarray}\label{closedform_neighbors0}
 && \Pi_{s_1, ..., s_n}(x_{s_{1}}, ..., x_{s_{n}})  \\
 &=& \sum_{i=2}^n \biggl[\biggl(\frac{\tilde{\pi}_{s_i}(x_{s_i})(x_{s_i}-\mu_{\tilde{\pi}_{s_i}})}{\sigma^2_{\tilde{\pi}_{s_i}}}\sum_{j=1}^{i-1} \prod_{k=1,k\not=j}^{i-1}\hat{\pi}_{s_k}(x_{s_k}) \times \frac{\beta_{s_i,s_j}\tilde{\pi}_{s_j}(x_{s_j})(x_{s_j}-\mu_{\tilde{\pi}_{s_j}})}{\sigma^2_{\tilde{\pi}_{s_j}}}\biggl)\times \prod_{k=i+1}^n \pi_{s_k}(x_{s_k}) \biggl] + \prod_{i=1}^n \pi_{s_i}(x_{s_i}) \nonumber
\end{eqnarray}
where \(\{\tilde{\pi}_{s_i}\}\) and \(\{\hat{\pi}_{s_i}\}\) are two auxiliary collections of pmfs, and we study conditions on \(\tilde{\pi}_{s_i}\), \(\hat{\pi}_{s_i}\), \(\pi_{s_i}\) and \(\beta_{s_i,s_j}\) that ensure the same random field \(X_S\) is constructed, regardless of site ordering. In addition, if \(\hat{\pi}_{s_i} \equiv \pi_{s_i}\) \(\forall 1\le i\le n\), then each pair of \textbf{uncorrelated} sites are actually \textbf{independent}. These results are given in Corllary \ref{cor:permutation:necesuffi}, Proposition \ref{permutation:necesuffm} and Corllary \ref{independenceofuncorrelatedneighbors}. The predominance followed from these results is to simulate true \textbf{Markov random fields} rather than just \textbf{correlated random fields} on a general site space \((S,\partial)\). For a given \((S,\partial)\), its neighborhood system \(\partial\) can be extended to \(\partial'\) (under \(\partial'\) all sites in \(S\) are neighbors of each other) by setting those pair of sites which are not neighbors of each other have covariances \textbf{0}, i.e., for each \(s\in S\), if any \(t\notin \partial(s)\), let \(\beta_{s,t}=0\). Herein, we develop a direct and flexible algorithm from Proposition \ref{permutation:necesuffm} to generate site-order invariant Markov random fields, i.e., the simulated Markov random field does not depend on the particular site order generating it. This subclass of random fields require stronger assumptions and have nicer properties, compared to the correlated random fields generated by Proposition \ref{mainprop:alg}.

The remainder of this note is laid out as follows: Section \ref{se:Notation_Background} contains our notation and the statement of
our main results, Proposition \ref{mainprop:alg}. Next, we explain our simulation algorithm in Section
\ref{se:Algor_Sim_Ran_Field}. In section \ref{se:Appl_Novel_Algorithm_to_Image_Analysis}, we summarize our prior applications of our simulation algorithm to image analysis. In Section \ref{se:AlgrthmConstrnts}, we explore its properties and constraints. We give the necessary and sufficient conditions of the regularity type for singleton case of \(A_{s_i}\) in Proposition \ref{prop:condition}. We also give a necessary and sufficient condition for the permutation property of a random field in Proposition \ref{permutation:necesuffi} and \ref{permutation:necesuffm}. An effective algorithm for generating Markov random field follows from these propositions. Section \ref{se:Proofs_Lemma_Proposition}, the Appendix, contains our proofs of Lemma \ref{ConnectedNS} and of Proposition \ref{mainprop:alg}.

\section{Notation and Background}\label{se:Notation_Background}

In this section, we recall some random field notation from Winkler \cite{Winkler},
introduce new concepts, and state a proposition from which our novel random fields and algorithm follow.
Let $S$ be a finite index set of \textbf{sites}; and for each site $s\in S$, $\mathbf{X}_s$ be a finite
\textbf{space of states} at site $s$.
For nonempty $A\subset S$, denote the space of configurations $x_A =(x_s)_{s\in A}$
on $A$ by $\mathcal{X}_A=\prod_{s\in A} \mathbf{X}_s$. If $A=S$, we abbreviate
$\mathcal{X}_S$ by $\mathcal{X}$, i.e., $\mathcal{X}=\prod_{s\in S} \mathbf{X}_s$.

Let $\Pi$ denote a \textbf{probability measure} or \textbf{distribution} on $\mathcal{X}$. If for every $x\in \mathcal{X}$, $\Pi(x)>0$, i.e., $\Pi$ is a strictly positive probability measure on $\mathcal{X}$, then $\Pi$ is called a \textbf{random field}. We also call the random vector $X=(X_s)_{s\in S}$ on the probability space $(\mathcal{X}, \Pi)$ a random field. For a nonempty $A\subset S$, define the projection map from $\mathcal{X}$ onto $\mathcal{X}_A$ as follows:
$$
  X_A: x \to x_A,
$$
where $x\in \mathcal{X}$ and $x_A \in \mathcal{X}_A$.

The definitions of neighborhood system and Markov random field follow directly from Winkler \cite{Winkler}:

{\defn\label{neighborhoodsytem} A \textbf{neighborhood system} $\partial =\{\partial(s): s\in S\}$ of $S$ is any
collection of subsets of $S$ that satisfies the conditions:
(i)  $s\notin \partial(s)$ for every $s\in S$ and (ii) $s\in \partial(t)$ if and only if $t\in \partial(s)$.
The sites $t\in \partial(s)$ are called \textbf{neighbors} of $s$.
}

{\defn\label{defnRandomField} The random field $\Pi$ is a \textbf{Markov random field} with respect to the neighborhood system \(\partial\) if for all \(x\in \mathcal{X}\),
\begin{equation*}
  \Pi(X_s = x_s | X_t = x_t, t\not= s)=\Pi(X_s = x_s| X_t = x_t, t\in \partial(s)).
\end{equation*}
}

We will introduce the notions of \textbf{base set} and \textbf{one pass simulatable} respectively. These concepts rely on the following definitions and discussion.
First, the neighborhood of nonempty set generalizes the neighborhood of single site:
{\defn\label{neighborhoodofA} The \textbf{neighborhood} of nonempty $A\subset S$ is $\displaystyle{\partial(A)=\biggl(\bigcup_{s\in A} \partial(s)\biggl)\setminus A}$, that is the neighbors of the sites in $A$ that are not part of $A$ themselves. For convenience, let \(\partial(\emptyset)=S\), i.e., the neighborhood of empty set \(\emptyset\) is all sites \(S\).
}

Next, we define \textbf{exclusiveness}, \textbf{connectedness} and \textbf{separateness} for subsets of site space \(S\) based on the neighborhood of nonempty set:
{\defn\label{exclusiveness,connectedness,separateness} Two subsets \(B_1, B_2\subset S\) are \textbf{exclusive} from each other if \([B_1\bigcup\partial(B_1)]\bigcap B_2=\emptyset\) (or \(B_1\bigcap [B_2\bigcup\partial(B_2)]=\emptyset\)). A subset \(B\subseteq S\) is \textbf{connected}, if for any nonempty proper subset \(A\subset B\), \(\partial(A)\bigcap B\not=\emptyset\) (If a subset \(B\subseteq S\) contains only one site or is empty, \(B\) is connected, since \(\partial(\emptyset) = S\)). A subset \(B\subseteq S\) is \textbf{separated} if it is not \textbf{connected}.
}

In the above definition, \([B_1\bigcup\partial(B_1)]\bigcap B_2=\emptyset\) and \(B_1\bigcap [B_2\bigcup\partial(B_2)]=\emptyset\) are equivalent, i.e., one implies the other. \textbf{Separateness} of a nonempty subset \(B\) implies that there exist two subsets \(B_1,B_2\) which are exclusive such that \(B = B_1\bigcup B_2\).

In some applications, e.g. image restoration and shape detection from shading, a site space
\(S\) is divided into two parts: good part and bad part, or known part and unknown part.
Herein, we let \( H \subset S\) denote the bad or unknown part, and \(H^C\circeq S\setminus H\) the good or known part of the site space \(S\).
For a nonempty subset \(H\) of a connected space $(S, \partial)$, we can order its sites
sequentially and associate each site with a connected set, which is a subset of the
neighborhood of the site:
{\lem\label{ConnectedNS} Assume that \((S,\partial)\) is a connected space of $N>2$ sites and \(H\subseteq S\) a subset of \(n\ge 1\) sites. Then, the sites in \(H\) can be ordered as a sequence \(\{s_i\}_{i=1}^n\) such that \(s_i\in
\partial(H^C\bigcup\{s_1,...,s_{i-1}\})\) for \(1\le i\le n\).
In addition, there exist unique \(m_i\ge 1\) and connected subsets \(\{B_{s_i}^j\}_{j=1}^{m_i}\) (called connected components) such that \(\displaystyle \partial(s_i)\bigcap [H^C\bigcup\{s_1,..., s_{i-1}\}]=\bigcup_{j=1}^{m_i} B_{s_i}^j\) and \(\{B_{s_i}^j\}_{j=1}^{m_i}\) are exclusive from each other. We choose one component among \(\{B_{s_i}^j\}_{j=1}^{m_i}\) and denote it by \(A_{s_i}\) (except \( A_{s_1}=\emptyset\) if \(H=S\)) for \(1\le i\le n\).
}

Lemma \ref{ConnectedNS} is proved in Section \ref{se:Proofs_Lemma_Proposition}.
Neither the site order nor the connected subsets \(\{A_{s_i}\}_{i=1}^n\) are
unique. Rather, we just assume henceforth that a particular setup has been chosen.

{\defn\label{validsetup} Suppose that \((S,\partial)\) is a connected space of $N>2$ sites
and \(H\subseteq S\) a subset of \(1\le n\le N\) sites.
Then, a \textbf{valid setup} for \(H\) is an ordered collection \(\{(s_i, A_{s_i})\}_{i=1}^n\)
such that $H=\bigcup_{i=1}^n \{s_i\}$, \(s_i\in
\partial(H^C\bigcup\{s_1,...,s_{i-1}\})\) for \(1\le i\le n\) and \(A_{s_i}\subseteq \partial (s_i)\bigcap [ H^C\bigcup\{s_1,..., s_{i-1}\}]\) is one connected component of \(\partial(s_i)\bigcap[H^C\bigcup\{s_1,...,s_{i-1}\}]\) for \(1\le i\le n\). The set $A_{s_i}$ is called the \textbf{base set} for site $s_i$.
}

{\examp When \(H=S\), a \textbf{valid setup} for \(H\) is an ordered
collection \(\{(s_i, A_{s_i})\}_{i=1}^N\) such that \(A_{s_1}=\emptyset\),
$s_{i}\in \partial(\{s_1, ..., s_{i-1}\})$ and $A_{s_i}$ being a connected component of $\partial (s_i)\bigcap \{s_1,..., s_{i-1}\}$ for $2\le i\le N$.
}

We are really interested in using a valid setup to simulate the unknown portion
of a random field.

{\defn The unknown states $H$ of a discrete random field $\Pi$ are
\label{SequentiallySimulatable}
\emph{one pass} or \emph{sequentially simulatable}
with valid setup \(\{(s_i, A_{s_i})\}_{i=1}^n\) if
\[
\Pi(X_{s_i}=x_{s_i}|X_{s_{i-1}}=x_{s_{i-1}},...,X_{s_{1}}=x_{s_{1}},X_{H^C}=x_{H^C})
  = \Pi( X_{s_i}=x_{s_i}|X_{A_{s_i}}=x_{A_{s_i}}).
\]
}
Then, by the multiplication rule one has that for any
configuration \(x\in \mathcal{X}\) is given:
\begin{eqnarray}\label{RandomFieldMultiplicationRule}
  \Pi(X_H=x_H|X_{H^C}=x_{H^C}) &=& \prod_{i=1}^n \Pi( X_{s_i}=x_{s_i}|X_{s_{i-1}}=x_{s_{i-1}},
  ...,X_{s_{1}}=x_{s_{1}},X_{H^C}=x_{H^C} ) \nonumber \\
 &=& \prod_{i=1}^n \Pi( X_{s_i}=x_{s_i}|X_{A_{s_i}}=x_{A_{s_i}}).
\end{eqnarray}
From this configuration probability, we can find the probability of any set of sites \(B\).
When \(B\subset\{s_1, ..., s_j\}\) for some \(j\), we need only compute the product of the first \(j\) terms
\begin{equation}
 \label{ComputeMarginal} \Pi(X_{B}=x_{B}|X_{H^C}=x_{H^C}) = \sum_{ s_k\notin B}\prod_{k=1}^{j} \Pi( X_{s_k}=x_{s_k}|X_{A_{s_k}}=x_{A_{s_k}}),
\end{equation}
and this formula can be used to compute
\(\Pi(X_{A_{s_i}}=x_{A_{s_i}}|X_{H^C}=x_{H^C})\) in (\ref{CondProb}) below.

Kolmogorov's consistency conditions consist of permutation condition and marginality condition. These conditions were applied to Kolmogorov's extension theorem to gurantee the unique existence (almost surely) of stochastic process that can degenerate to given families of random vectors. These conditions also light up those good properties that Markov random fields satisfy. For each integer \(1<n\le N\), let \(M_n\) contain \(n\) elements and \(G_{M_n}\) be the symmetric group consisting of all permutations on \(M_n\). Then we define the consistency conditions for discrete random field as follows:
{\defn\label{ConsistencyCondi} Assume that \((S,\partial)\) is a connected space of \(N>2\) sites, \(H\)
is a subset of \(n\ge 1\) sites with valid setup \(\{(s_i, A_{s_i})\}_{i=1}^n\)
and $X_{H^C}$ be the known portion of the random field. For each \(1\le k\le n\), and each \(M_k=\{i_1,...,i_k\}\subseteq \{1,\cdots,n\}\) and \(g\in G_{M_k}\), let \(\Pi_{s_{g(i_1)},...,s_{g(i_k)}}(\cdot|X_{H^C}=x_{H^C})\) (if it can be defined the same way as \(\Pi_{s_{i_1},...,s_{i_k}}(\cdot|X_{H^C}=x_{H^C})\)) be a probability mass function on \(\mathbf{X}_{s_{g(i_1)},..., s_{g(i_k)}}=\prod_{u=1}^k\mathbf{X}_{s_{g(i_u)}}\). The two consistency conditions for the collection of probability mass functions \(\{\Pi_{s_{g(i_1)},...,s_{g(i_k)}}(\cdot|X_{H^C}=x_{H^C}): g\in G_{M_k}, M_k=\{i_1,...,i_k\}\subseteq \{1,\cdots,n\}, 1\le k\le n\}\) are listed as follows.
\begin{enumerate}
\item \textbf{permutation}: for each permutation \(g\in G_{M_k}\), \(M_k=\{i_1,...,i_k\}\subseteq \{1,\cdots,n\}\) (\(1\le  k\le n\)), and \(x_{s_{i_u}}\in \mathbf{X}_{s_{i_u}}\) (\(1\le u\le k\)),
\begin{equation*}
    \Pi_{s_{g(i_1)},...,s_{g(i_k)}} (x_{s_{g(i_1)}}, ...,x_{s_{g(i_k)}}|X_{H^C}=x_{H^C}) = \Pi_{s_{i_1},...,s_{i_k}}(x_{s_{i_1}}, ...,x_{s_{i_k}}|X_{H^C}=x_{H^C});
\end{equation*}
\item \textbf{marginality}: for each permutation \(g\in G_{M_k}\), \(M_k=\{i_1,...,i_k\}\subseteq \{1,\cdots,n\}\) (\(2\le  k\le n\)) and \(1\le j\le k\),
\begin{eqnarray*}
    &&\sum_{x_{s_{g(i_j)}}\in \mathbf{X}_{s_{g(i_j)}}}\Pi_{s_{g(i_1)},...,s_{g(i_k)}}(x_{s_{g(i_1)}}, ...,x_{s_{g(i_k)}}|X_{H^C}=x_{H^C}) \\
    &=& \Pi_{s_{g(i_1)},...,s_{g(i_{j-1})},s_{g(i_{j+1})},...,s_{g(i_k)}}(x_{s_{g(i_1)}}, ..., x_{s_{g(i_{j-1})}}, x_{s_{g(i_{j+1})}},..., x_{s_{g(i_k)}}|X_{H^C}=x_{H^C}),
\end{eqnarray*}
where \(x_{s_{g(i_u)}}\in \mathbf{X}_{s_{g(i_u)}}\) (\(1\le u\not=j\le k\)).
\end{enumerate}
 }

Note that in the above definition, it can happen that some \(\Pi_{s_{g(i_1)},...,s_{g(i_k)}}(\cdot|X_{H^C}=x_{H^C})\) are not defined. For example, for \(g\) on \(M_n=\{1,\cdots,n\}\), the sites in \(H\) can be ordered by \(g\) as \(\{s_{g(i)}\}_{i=1}^n\), but it may fail that \(s_{g(i)}\in \partial(H^C\bigcup\{s_{g(1)},\cdots, s_{g(i-1)}\})\) for some \(1\le i\le n\). In such case, we do not have \(\Pi_{s_{g(1)},...,s_{g(n)}}\) defined. But if \(\Pi_{s_{g(1)},...,s_{g(n)}}\) is defined indeed, then the permutation property states that the random fields on \(H\) generated by the order \(X_{s_1}\), ..., \(X_{s_n}\) and the order \(X_{s_{g(1)}}\), ..., \(X_{s_{g(n)}}\) are the same. So the condition for permutation property assures that \(X_H\) on \(H\) can be simulated in any site order, given \(X_H\) can be generated in such order. Note also permutation property is not required or recommended for some type of random fields such as discrete time series. The marginality condition ensures that random field  \(X_{H_1}\) generated directly on a proper subset \(H_1\subset H\) is the same as that degenerated from \(X_H\).

When simulating \(X_{s_i}\), one has access to  \(X_{H^C}\bigcup
\{X_{s_1}, ..., X_{s_{i-1}}\}\), the known sites and the sites already simulated,
but need only use \(X_{A_{s_i}}\).
This reduction is the key that makes one-pass sequential simulation effective
when the base sets are relatively small compared to \(H\).
We will explain how this can be done when one is just interested in
simulating a field with given marginal probabilities and certain covariances. For the case \(H=S\), since \(H^C=\emptyset\), the conditioning \(X_{H^C}=x_{H^C}\) disappears from (\ref{RandomFieldMultiplicationRule})
 and (\ref{ComputeMarginal}), leading to an easier-to-assimulate means of constructing a random field. Conversely, the case \(H\subset S\)
 is advantageous for real applications.

Herein, we simulate random fields with given marginal probabilities for
sites and given covariances between nearby (meaning within the base sets however they are defined)
sites for a subset \(H\) of a connected space $(S, \partial)$ with $N$ sites.
We assume a valid setup \(\{(s_i, A_{s_i})\}_{i=1}^n\) for \(H\). Our algorithm constructs $X_H$
with the given marginal probabilities
\(\{\pi_{s_i}(x_{s_i}): x_{s_i}\in \mathbf{X}_{s_i}\}_{i=1}^n\)
and the given covariances between nearby sites
$\{\beta_{s_i,t_i}: t_i\in A_{s_i} \}$ for $1\le i\le n$.
(It is assumed a priori that these marginal conditions hold
within the used set of known sites $(\cup_{i=1}^n A_{s_i})\setminus H$. It is also assumed that \(\beta_{s_i,t_i}=\beta_{t_i,s_i}\) for \(t_i\in A_{s_i}\) (\(1\le i\le n\)) since \(\beta_{s_i,t_i}\) will denote covariance between \(X_{s_i}\) and \(X_{t_i}\)). We assign conditional probabilities
\(\Pi(X_{s_i}=x_{s_i}|\text{ } X_{A_{s_i}}=x_{A_{s_i}})\),
such that we maintain the desired covariances and marginal probabilities
as we include the unknown sites:
\begin{prop}
\label{mainprop:alg}
Assume that \((S,\partial)\) is a connected space of \(N>2\) sites, \(H\)
is a subset of \(n\ge 1\) sites with valid setup \(\{(s_i, A_{s_i})\}_{i=1}^n\)
and $X_{H^C}$ is the known portion of the random field.
Suppose further that
\(
\displaystyle \{\tilde{\pi}_{s}(x_{s}): x_{s}\in \mathbf{X}_{s}, s \in S\}
\)
and
\(
\displaystyle \{\hat{\pi}_{s}(x_{s}): x_{s}\in \mathbf{X}_{s}, s\in S\}
\)
are two sets of pmfs. Assume that
\(
\displaystyle \{\pi_{s}(x_{s}): x_{s}\in \mathbf{X}_{s}, s\in S\}
\)
are positive pmfs
and \(\{\beta_{s_i,t_i}: t_i\in A_{s_i}, 1\le i\le n\}\) are numbers such that
the right hand side (RHS) of (\ref{CondProb}) is in [0,1] for all $i$.
Form the conditional probabilities starting with $i=1$ recursively as
\begin{eqnarray}\label{CondProb}
&&\Pi(\left. X_{s_i}=x_{s_i}\right| X_{A_{s_i}}=x_{A_{s_i}})  \\
&=& \pi_{s_i}(x_{s_i})+\frac{\tilde{\pi}_{s_i}(x_{s_i})(x_{s_i}-\mu_{\tilde{\pi}_{s_i}})}{\sigma^2_{\tilde{\pi}_{s_i}}} \sum_{t_i\in A_{s_i}}\biggl(\prod_{u_i\in A_{s_i}\setminus \{t_i\}}\hat{\pi}_{u_i}(x_{u_i})\biggl)\cdot \frac{\displaystyle{
\beta_{s_i,t_i}\tilde{\pi}_{t_i}(x_{t_i})(x_{t_i}-\mu_{\tilde{\pi}_{t_i}})}}
{\displaystyle{\sigma^2_{\tilde{\pi}_{t_i}}}
\Pi(X_{A_{s_i}}=x_{A_{s_i}}| X_{H^C}=x_{H^C})} \nonumber
\end{eqnarray}
for each \(x_{s_i}\in \mathbf{X}_{s_i}\) and
\(x_{A_{s_i}} \in \mathbf{X}_{A_{s_i}}\) (\(1\le i\le n\)),
where $\mu_{\tilde{\pi}_{s}}=\displaystyle\sum_{x_{s}\in \mathbf{X}_{s}} \tilde{\pi}_{s}(x_s) x_s$ and
$\sigma^2_{\tilde{\pi}_{s}}=\displaystyle \sum_{x_{s}\in \mathbf{X}_{s}} \tilde{\pi}_{s}(x_s)(x_s-\mu_{\tilde{\pi}_{s}})^2$ (\(s\in S\)).
Then, there is a probability measure $\Pi$ on
$\mathcal{X}_{(\cup_{i=1}^n A_{s_i})\cup H}$ consistent with
(\ref{CondProb}) that has marginal probabilities \(\{\pi_{s_i}\}\)
and covariances
\(\text{cov}(X_{s_i},X_{t_i})=\beta_{s_i,t_i}\) for all $t_i\in A_{s_i}, 1\le i\le n$.
\end{prop}

{\rem\label{HisSubsetofS} Proposition \ref{mainprop:alg} can be used in image smoothing when \(H\subset S\).
For example, an image \(X\) can be smoothed in a way as follows:
let \(H\) be the set of pixels with states consisting of ``sharp" or
undesirable values, and replace those values by simulated ones, using
Proposition \ref{mainprop:alg}.
$H^C$ is the portion of the picture that does not require
smoothing.}

{\rem The special case \(H=S\) of Proposition \ref{mainprop:alg} deserves particular attention. When \(H=S\), \(A_{s_1}=\emptyset\), the second term on the right hand side of (\ref{CondProb}) disappears since the summation is over zero terms. \(\Pi(X_{A_{s_i}}=x_{A_{s_i}}| X_{H^C}=x_{H^C})\) in (\ref{CondProb}) should be replaced by \(\Pi(X_{A_{s_i}}=x_{A_{s_i}})\) because \(H^C=S^C=\emptyset\). Different from the case \(H\subset S\) in Remark \ref{HisSubsetofS}, Proposition \ref{mainprop:alg} with \(H=S\) is used to generate random field $X$ on the site space $S$: \(\Pi(X=x) = \displaystyle\prod_{i=1}^N \Pi( X_{s_i}=x_{s_i}|X_{A_{s_i}}=x_{A_{s_i}})\).
}

{\rem In Proposition \ref{mainprop:alg}, we assumed that \( \pi_{s}(x_{s})>0: \forall x_{s}\in \mathbf{X}_{s}\) for each \(s\in S\). Note that \(\mathbf{X}_{s}\) can be different for each \(s\in S\). For given \(s\), if there exists a \(x_{s}\in \mathbf{X}_{s}\) such that \(\pi_s(x_s)=0\), we may deem it uninteresting and replace \(\mathbf{X}_{s}\) with \(\mathbf{X}_{s}\setminus \{x_s\}\). Therefore the positive probability mass function assumption of \(\pi_s\) is also a convention. But for \(\tilde{\pi}_{s}\) and \(\hat{\pi}_{s}\), we do not need to have this assumption, i.e. the support of \(\tilde{\pi}_{s}(x_{s})\) or \(\hat{\pi}_{s}\) can be a proper subset of \(\mathbf{X}_{s}\).
}

{\rem In general we will match covariances between a site \(s\) and sites in \(A_s^{cov}=A_s\cup  \{t:s\in A_t\}\). This is one of the connected components in \(\partial(s)\).
}

{\rem\label{secondtherm} By definition, \textbf{connectedness} is based upon the notion of neighbor, and so is the Markov property of a Markov random field. The connectedness condition of the site space \(S\) guarantees that we can use the accumulated information about the neighbor states \(X_{t_i}\), \(t_i\in A_{s_i}\), to generate the state \(X_{s_i}\).
}

A random field generated by Proposition \ref{mainprop:alg} is a correlated random field. Indeed, one value of this proposition is the assertion that there are one-pass simulatable correlated random fields that match a given collection of marginal probabilities and covariances. What's more, we will state conditions in Section \ref{se:AlgrthmConstrnts} that convert a correlated random field into a Markov random field. We call these one-pass simulatable correlated random fields generated by Proposition \ref{mainprop:alg} the \textbf{KNW Random Fields} for ease of future reference.

{\rem Proposition \ref{mainprop:alg} has a few interesting and important special cases:
\begin{enumerate}
\item If \(\hat{\pi}_{s}(x_{s})=\tilde{\pi}_{s}(x_{s}), x_{s}\in \mathbf{X}_{s}\) for all \(s\in S\), then (\ref{CondProb}) becomes
\begin{equation}\label{CondProbEqualSets}
\Pi(\left. X_{s_i}=x_{s_i}\right| X_{A_{s_i}}=x_{A_{s_i}})
= \pi_{s_i}(x_{s_i})+\bigg(\prod_{u_i\in \bar{A}_{s_i}}\tilde{\pi}_{u_i}(x_{u_i})\biggl)\cdot\sum_{t_i\in A_{s_i}}\frac{(x_{s_i}-\mu_{\tilde{\pi}_{s_i}})
\beta_{s_i,t_i}(x_{t_i}-\mu_{\tilde{\pi}_{t_i}})}{\sigma^2_{\tilde{\pi}_{s_i}}\sigma^2_{\tilde{\pi}_{t_i}}\Pi(X_{A_{s_i}}=x_{A_{s_i}}| X_{H^C}=x_{H^C})}
\end{equation}
where \(\bar{A}_{s_i}=A_{s_i}\bigcup\{s_i\}\). In (\ref{CondProbEqualSets}), two auxiliary collections of pmfs are reduced to one.

\item If \(\hat{\pi}_{s}(x_{s})=\tilde{\pi}_{s}(x_{s})=\pi_{s}(x_{s}), x_{s}\in \mathbf{X}_{s}\) for all \(s\in S\), then (\ref{CondProb}) is
\begin{eqnarray}\label{CondProbNormalized}
&&\Pi(\left. X_{s_i}=x_{s_i}\right| X_{A_{s_i}}=x_{A_{s_i}}) \nonumber \\
&=& \pi_{s_i}(x_{s_i})\biggl[1+\bigg(\prod_{u_i\in A_{s_i}}\pi_{u_i}(x_{u_i})\biggl)\cdot\sum_{t_i\in A_{s_i}}\frac{(x_{s_i}-\mu_{\pi_{s_i}})
\beta_{s_i,t_i}(x_{t_i}-\mu_{\pi_{t_i}})}{\sigma^2_{\pi_{s_i}}\sigma^2_{\pi_{t_i}}\Pi(X_{A_{s_i}}=x_{A_{s_i}}| X_{H^C}=x_{H^C})}\biggl]
\end{eqnarray}
There are no auxiliary collections of pmfs in (\ref{CondProbNormalized}).

\item If \(\hat{\pi}_{s}(x_{s})=\tilde{\pi}_{s}(x_{s})=\frac{1}{d_s}, x_{s}\in \mathbf{X}_{s}\) for all \(s\in S\) where \(d_s\) is the cardinality of \(\mathbf{X}_{s}\), then (\ref{CondProb}) has the following form
\begin{equation}\label{CondProbUniform}
\Pi(\left. X_{s_i}=x_{s_i}\right| X_{A_{s_i}}=x_{A_{s_i}})
= \pi_{s_i}(x_{s_i})+\sum_{t_i\in A_{s_i}}\frac{(x_{s_i}-\bar{\mu}_{s_i})
\beta_{s_i,t_i}(x_{t_i}-\bar{\mu}_{t_i})}{D_{\bar{A}_{s_i}}\bar{\sigma}^2_{s_i}\bar{\sigma}^2_{t_i}\Pi(X_{A_{s_i}}=x_{A_{s_i}}| X_{H^C}=x_{H^C})}
\end{equation}
where \(\displaystyle\bar{\mu}_s=\frac{1}{d_s}\sum_{x_{s}\in \mathbf{X}_{s}}x_s\) and \(\displaystyle\bar{\sigma}^2_s=\frac{1}{d_s}\sum_{x_{s}\in \mathbf{X}_{s}}(x_s-\bar{\mu}_s)^2\) for all \(s\in S\) and \(\displaystyle D_{\bar{A}_{s_i}}=\prod_{u_i\in \bar{A}_{s_i}}d_{u_i}\). Here two collections of auxiliary pmfs take same discrete uniform pmfs respectively. The simplicity in (\ref{CondProbUniform}) reduces the computation of conditional probabilities. Notice the \(\bar{\mu}_s\) and \(\bar{\sigma}^2_s\) calculations are simplified.

\item If \(\hat{\pi}_{s}(x_{s})=\frac{1}{d_s}\) and \(\tilde{\pi}_{s}(x_{s})=\pi_{s}(x_{s})\), \(\forall x_{s}\in \mathbf{X}_{s}\) for all \(s\in S\) with \(d_s\) being the cardinality of \(\mathbf{X}_{s}\), then (\ref{CondProb}) is changed to
\begin{equation}\label{CondProbUniform2}
\Pi(\left. X_{s_i}=x_{s_i}\right| X_{A_{s_i}}=x_{A_{s_i}})
= \pi_{s_i}(x_{s_i})\biggl[1+\sum_{t_i\in A_{s_i}}\frac{(x_{s_i}-\mu_{\pi_{s_i}})\beta_{s_i,t_i}\pi_{t_i}(x_{t_i})(x_{t_i}-\mu_{\pi_{t_i}})}{D_{A_{s_i}\setminus\{t_i\}}\sigma^2_{\pi_{s_i}}\sigma^2_{\pi_{t_i}}\Pi(X_{A_{s_i}}=x_{A_{s_i}}| X_{H^C}=x_{H^C})}\biggl]
\end{equation}
where \(\displaystyle D_{A_{s_i}\setminus\{t_i\}}=\prod_{u_i\in A_{s_i}\setminus\{t_i\}}d_{u_i}\). Here one collection of auxiliary pmfs take discrete uniform pmfs and another collection is identitical to the prescribed \(\{\pi_s\}\).

\item If we assume the sufficient condition for permutation property of a random field (see Corollary \ref{cor:permutation:necesuffi}): \(\hat{\pi}_{s}(x_{s})=\pi_{s}(x_{s}), \forall x_{s}\in \mathbf{X}_{s}\) for all \(s\in S\), then (\ref{CondProb}) becomes
\begin{eqnarray}\label{CondPermutationProperty}
&&\Pi(\left. X_{s_i}=x_{s_i}\right| X_{A_{s_i}}=x_{A_{s_i}}) \nonumber \\
&=& \pi_{s_i}(x_{s_i})+\frac{\tilde{\pi}_{s_i}(x_{s_i})(x_{s_i}-\mu_{\tilde{\pi}_{s_i}})}{\sigma^2_{\tilde{\pi}_{s_i}}} \sum_{t_i\in A_{s_i}}\biggl(\prod_{u_i\in A_{s_i}\setminus \{t_i\}}\pi_{u_i}(x_{u_i})\biggl)\cdot \frac{\displaystyle{
\beta_{s_i,t_i}\tilde{\pi}_{t_i}(x_{t_i})(x_{t_i}-\mu_{\tilde{\pi}_{t_i}})}}
{\displaystyle{\sigma^2_{\tilde{\pi}_{t_i}}}
\Pi(X_{A_{s_i}}=x_{A_{s_i}}| X_{H^C}=x_{H^C})}
\end{eqnarray}

\item Assume that the space of states \(\mathbf{X}_s\) for all \(s\in S\) is same and is denoted by \(\mathbf{X}\). If \(\hat{\pi}_{s}(x_{s})=\hat{\pi}(x_{s})\) and \(\tilde{\pi}_{s}(x_{s})=\tilde{\pi}(x_{s})\), \(x_{s}\in \mathbf{X}\) for all sites \(s\in S\), then (\ref{CondProb}) is adapted to
\begin{eqnarray}\label{CondProbTwoModifiers}
&&\Pi(\left. X_{s_i}=x_{s_i}\right| X_{A_{s_i}}=x_{A_{s_i}}) \nonumber \\
&=& \pi_{s_i}(x_{s_i})+\frac{\tilde{\pi}(x_{s_i})(x_{s_i}-\mu_{\tilde{\pi}})}{\sigma^2_{\tilde{\pi}}} \sum_{t_i\in A_{s_i}}\biggl(\prod_{{u_i}\in A_{s_i}\setminus \{t_i\}}\hat{\pi}(x_{u_i})\biggl)\cdot \frac{\displaystyle{
\beta_{s_i,t_i}\tilde{\pi}(x_{t_i})(x_{t_i}-\mu_{\tilde{\pi}})}}
{\displaystyle{\sigma^2_{\tilde{\pi}}}
\Pi(X_{A_{s_i}}=x_{A_{s_i}}| X_{H^C}=x_{H^C})}
\end{eqnarray}
where $\mu_{\tilde{\pi}}=\displaystyle\sum_{x_{s}\in \mathbf{X}} \tilde{\pi}(x_s) x_s$ and
$\sigma^2_{\tilde{\pi}}=\displaystyle \sum_{x_{s}\in \mathbf{X}} \tilde{\pi}(x_s)(x_s-\mu_{\tilde{\pi}})^2$ (\(s\in S\)). In this case, auxiliary pmfs \(\{\hat{\pi}_s\}\) are identically distributed, so are \(\{\tilde{\pi}_s\}\).

\item If we combine the assumptions for formula (\ref{CondProbUniform}) and (\ref{CondProbTwoModifiers}) together, i.e., \(\mathbf{X}_s=\mathbf{X}\) and \(\hat{\pi}_{s}(x_{s})=\tilde{\pi}_{s}(x_{s})=\frac{1}{d}, x_{s}\in \mathbf{X}\) for all \(s\in S\) where \(d\) is the cardinality of \(\mathbf{X}\), then (\ref{CondProb}) takes the following simple form
\begin{equation}\label{Captchaformula}
\Pi(\left. X_{s_i}=x_{s_i}\right| X_{A_{s_i}}=x_{A_{s_i}})
= \pi_{s_i}(x_{s_i})+\sum_{t_i\in A_{s_i}}\frac{(x_{s_i}-\bar{\mu})
\beta_{s_i,t_i}(x_{t_i}-\bar{\mu})}{d^{|A_{s_i}|+1}(\bar{\sigma}^2)^2\Pi(X_{A_{s_i}}=x_{A_{s_i}}| X_{H^C}=x_{H^C})}
\end{equation}
where \(\displaystyle\bar{\mu}=\frac{1}{d}\sum_{x_{s}\in \mathbf{X}}x_s\), \(\displaystyle\bar{\sigma}^2=\frac{1}{d}\sum_{x_{s}\in \mathbf{X}}(x_s-\bar{\mu})^2\) \((s\in S)\) and \(|A_{s_i}|\) is the cardinality of \(A_{s_i}\). Formula (\ref{Captchaformula}) is used in \cite{KoNeWu} for Captcha generation.
\end{enumerate}
}

Insomuch as the case $H\not=S$ only involves a notational change, the remainder of this note will only consider the case \(H=S\) and will not state it anymore. For example, in the appendix we will only prove Proposition \ref{mainprop:alg} for the case $H=S$.

{\examp \label{example1} Let \(S=\{1,2,3,4,5\}\), \(s_i=i\)
for \(1\le i \le 5\) and \(\partial\) contain the following subsets
of \(S\): \(\partial(1)=\{2,4\}\), \(\partial(2)=\{1,3,5\}\),
\(\partial(3)=\{2,4,5\}\), \(\partial(4)=\{1,3,5\}\) and
\(\partial(5)=\{2,3,4\}\).
Then, \(\partial\) satisfies the conditions of
a neighborhood system, i.e. (i) \(i\notin \partial(i)\) for \(i\in S\)
and (ii) \(i\in \partial(j)\) if and only if \(j\in \partial(i)\). The connected space \((S,\partial)\) is drawn as follows: each site (\(1\le i\le 5\)) is represented by a node and each pair of neighbors are connected by a dashed edge.

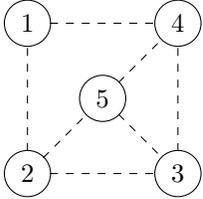
\begin{tikzpicture}
\path ( 0,2) node[shape=circle,draw] {1}
(0,0) node[shape=circle,draw] {2}
(2,0) node[shape=circle,draw] {3}
(2,2) node[shape=circle,draw] {4}
(1,1) node[shape=circle,draw] {5};
\draw [dashed] (0,0.3) -- (0,1.7);
\draw [dashed] (0.3,0) -- (1.7,0);
\draw [dashed] (2,0.3) -- (2,1.7);
\draw [dashed] (0.3,2) -- (1.7,2);
\draw [dashed] (0.21,0.21) -- (0.79,0.79);
\draw [dashed] (1.21,1.21) -- (1.79,1.79);
\draw [dashed] (1.21,0.79) -- (1.79,0.21);
\end{tikzpicture}
}

By previous discussion, we can take \(A_1=\emptyset\),
\(A_2=\partial(2)\bigcap\{1\}=\{1\}\),
\(A_3=\partial(3)\bigcap\{1,2\}=\{2\}\),
\(A_4\subset\partial(4)\bigcap\{1,2,3\}=\{1,3\}\),
\(A_5=\partial(5)\bigcap\{1,2,3,4\}=\{2,3,4\}\).
Now, \(\{1,3\}\) is not
connected because \(\{1\}\) is a proper subset and
\(\partial(\{1\})\bigcap \{1,3\}=\emptyset\).
Hence, we choose \(A_4 = \{3\}\)
(We could have chosen \(A_4 = \{1\}\) as well). Note that
\(A_5=\{2,3,4\}\) is a connected set with respect to the
neighborhood system \(\partial\).

Let the common space of states for each site be
\(\mathbf{X}=\{-1,0,1\}\). To construct a probability measure $\Pi$
on \( \mathcal{X}= \{-1,0,1\}^5 \), we use Proposition \ref{mainprop:alg} to compute the conditional probabilities: \(\Pi(X_1=x_1)\),
\(\Pi(X_2=x_2|X_1=x_1)\), \(\Pi(X_3=x_3|X_2=x_2)\),
\(\Pi(X_4=x_4|X_3=x_3)\), and \(\Pi(X_5=x_5|X_{A_5}=x_{A_5})\) where
\(x_i\in \{-1,0,1\}\) for \(1\le i\le 5\) and \( x_{A_5} \in
\{-1,0,1\}^3\). To compute \(\Pi(X_5=x_5|X_{A_5}=x_{A_5})\), we have to compute
\(\Pi(X_{A_5}=x_{A_5})\) as a prerequisite, using Proposition \ref{mainprop:alg}. By (\ref{RandomFieldMultiplicationRule}) and (\ref{ComputeMarginal}), first
we compute the joint probabilities
\begin{eqnarray*}
  \Pi(X_4=x_4, X_3=x_3, X_2=x_2, X_1=x_1)&=&\Pi(X_4=x_4|X_3=x_3)\Pi(X_3=x_3|X_2=x_2) \\
     && \times \Pi(X_2=x_2|X_1=x_1)\Pi(X_1=x_1),
\end{eqnarray*}
Since \(A_2=\{1\}\), \(A_3=\{2\}\) and \(A_4=\{3\}\), we get
\[
  \Pi(X_{A_5}=x_{A_5})=\sum_{x_1\in \{-1,0,1\}}\Pi(X_4=x_4, X_3=x_3, X_2=x_2, X_1=x_1).
\]

Alternatively, with the connectedness of \(A_5\), we can treat \(A_5\) as \(S\), use Proposition
\ref{mainprop:alg} again, and compute \(\Pi(X_{A_5}=x_{A_5})\) as
follows:
\[
  \Pi(X_{A_5}=x_{A_5})=\Pi(X_{4}=x_{4}|X_{3}=x_{3})\Pi(X_{3}=x_{3}|X_{2}=x_{2})\Pi(X_{2}=x_{2}).
\]
This later alternative method is more efficient than the former one. \qed

{\rem\label{prop} Example \ref{example1} shows how to use the
connectedness assumption on \(A_{s_i}\) (\( 2\le i\le N\)) in Proposition \ref{mainprop:alg}. By Definition \ref{validsetup}, $A_{s_i}$ ($2\le i\le N$) is chosen to be one of the connected components of $\partial(s_i)\bigcap \{s_1, ..., s_{i-1}\}$. For each site \(s_i\) (\(1\le i\le N\)), let \(A^{cov}_{s_i}=A_{s_i}\bigcup\{s_j: s_i\in A_{s_j}, 1\le j\le N\}\). Herein, \(A^{cov}_{s_i}\) is the largest subset of \(\partial(s_i)\) containing $A_{s_i}$ such that the covariances between \(X_{s_i}\) and \(X_{t_i}\) (\(t_i\in A^{cov}_{s_i}\)) are matched. For the choice \(A_i\) (\(1\le i\le 5\)) in Example \ref{example1}, \(A^{cov}_1 = \{2\}\), \(A^{cov}_2 = \{1,3,5\}\), \(A^{cov}_3 = \{2,4,5\}\), \(A^{cov}_4 = \{3,5\}\) and \(A^{cov}_5 = \{2,3,4\}\). In the following diagram, those neighbors with matched covariances are connected by thick edges, and the covariance between \(X_1\) and \(X_4\) are not matched:

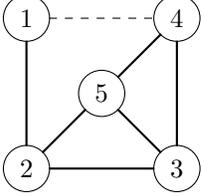
\begin{tikzpicture}
\path ( 0,2) node[shape=circle,draw] {1}
(0,0) node[shape=circle,draw] {2}
(2,0) node[shape=circle,draw] {3}
(2,2) node[shape=circle,draw] {4}
(1,1) node[shape=circle,draw] {5};
\draw [thick] (0,0.3) -- (0,1.7);
\draw [thick] (0.3,0) -- (1.7,0);
\draw [thick] (2,0.3) -- (2,1.7);
\draw [dashed] (0.3,2) -- (1.7,2);
\draw [thick] (0.21,0.21) -- (0.79,0.79);
\draw [thick] (1.21,1.21) -- (1.79,1.79);
\draw [thick] (1.21,0.79) -- (1.79,0.21);
\end{tikzpicture}
}

{\rem\label{Asiisonlycompnent} If \(S\) is connected and \(A_{s_i}=\partial
(s_i)\bigcap \{s_1,..., s_{i-1}\}\) for \( 2\le i\le N\), then \(A^{cov}_{s_i}=\partial
(s_i)\), \(\forall 1\le i\le N\). Here \(A_{s_i}\) is the only connected component of \(\partial
(s_i)\bigcap \{s_1,..., s_{i-1}\}\). Herein we can prescribe \(\{\beta_{s,t}: t\in \partial(s), s \in S\}\) such that the random field \(X\) constructed by
Proposition \ref{mainprop:alg} satisfies \(cov(X_s,X_t)=\beta_{s,t}\) for any pair of neighbors \(s\), \(t\in S\)
(without reference to the sets $A_s$).
}

As introduced in Section \ref{sec:Introduction}, random fields have typical applications in two-dimensional space (e.g. X-ray imaging) and three-dimensional space (e.g. human brain imaging and mammography). Herein, we present an example of two-dimensional site space and two types of its natural neighborhood systems. The three-dimensional analogue  can be formulated accrodingly.

{\examp\label{twodimensionalimagespace} Let \(S\) be the space of pixels on an image of size \(M\times N\), that is to say, \(S=\{(i,j):1\le i\le M, 1\le j\le N\}\), where \(M,N\in \mathbb{N}\). Recall the \textbf{square distance} \(\rho\) defined on \(\mathbb{R}^2\): \(\rho(P,Q)=\max(|x_2-x_1|, |y_2-y_1|)\) where \(P(x_1, y_1), Q(x_2,y_2)\in \mathbb{R}^2\). Fix \(\ell\in \mathbb{N}\), \(\ell\) is the radius of neighborhood system \(\partial_{\ell}\) on \(S\) defined as follows: for each \((i,j)\in S\),
\[
\partial_{\ell}((i,j))=\{(u,v)\in S: 0<\rho((i,j),(u,v))\le \ell\}.
\]
Therefore, \((S,\partial_{\ell})\) is a connected space. There are \(MN\) pixels in \(S\), and we enumerate these pixels in the following sequence: \(s_1=(1,1), s_2=(2,1), ..., s_M=(M,1)\), \(s_{M+1}=(1,2), s_{M+2}=(2,2), ..., s_{2M}=(M,2), ..., s_{(N-1)M+1}=(1,N)\), \(s_{(N-1)M+2}=(2,N), ..., s_{NM}=(M,N)\), i.e., we list the pixels on the first column of the image in increasing order, then the second column until the last column. The \((i,j)^{th}\) pixel on the image is the \(((j-1)M+i)^{th}\) element in the sequence, i.e.,
\( s_{(j-1)M+i}=(i,j) \).
For each \((i,j)\in S\), let
\begin{equation}\label{defnbiggestbaseset}
 A^{\ell}_{(i,j)}\triangleq\partial_{\ell}((i,j))\bigcap\{s_1, s_2, ..., s_{(j-1)M+i-1}\}
=\partial_{\ell}(s_{(j-1)M+i})\bigcap\{s_1, s_2, ..., s_{(j-1)M+i-1}\}
\end{equation}
be the points in \(S\) within the square distance \(\ell\) from \((i,j)\) that are either to the left or directly above of \((i,j)\). The restriction to \(\{s_1, s_2, ..., s_{(j-1)M+i-1}\}\) ensures that only the previously simulated pixels that are to the left or directly above of \((i,j)\) are used to generate the next pixel \((i,j)\).
One can verify that
\begin{equation}\label{biggestbaseset}
 A^{\ell}_{(i,j)}=\{(u,v):i-(i-1)\wedge\ell\le u\le (i+l)\wedge M, j-(j-1)\wedge\ell\le v\le j-1\}\bigcup\{(u,j):i-(i-1)\wedge\ell\le u<i\}.
\end{equation}
It can also be shown that \(\partial_{\ell}((i,j))\bigcap\{s_1, s_2, ..., s_{(j-1)M+i-1}\}\) contains only one connected component, i.e., $A^{\ell}_{(i,j)}$.

Example \ref{twodimensionalimagespace} shows that the abstract connected site space \((S,\partial)\) can have natural instances from real applications. The sequence of pixels enumerated above satisfies that $s_{k}\in \partial(\{s_1, ..., s_{k-1}\})$ for $2\le k\le MN$, meanwhile \(\{A^{\ell}_{(i,j)}\}_{(i,j)\in S\setminus\{(1,1)\}}\) satisfy that  \(A^{\ell}_{(i,j)}=\partial_{\ell}((i,j))\bigcap\{s_1, s_2, ..., s_{(j-1)M+i-1}\}\). Therefore, by Remark \ref{Asiisonlycompnent}, \(A^{\ell, cov}_{(i,j)}=\partial_{\ell}((i,j))\), \(\forall (i,j)\in S\).

Kouritzin et al \cite{KoNeWu} defined the same site space \(S\) as that in Example \ref{twodimensionalimagespace}, but replaced the square distance with \textbf{Euclidean distance} and \(\ell\in \mathbb{R}^+\) is not necessarily an integer. When \(\ell\ge \sqrt{2}\), the base set \(A^{\ell}_{(i,j)}\) can still be defined by (\ref{defnbiggestbaseset}) albeit (\ref{biggestbaseset}) does not hold. Given \(\ell\ge \sqrt{2}\), we still have \(A^{\ell, cov}_{(i,j)}=\partial_{\ell}((i,j))\), \(\forall (i,j)\in S\). Nevertheless, when \(1\le \ell < \sqrt{2}\), \(\partial_{\ell}((i,j))\bigcap\{s_1, s_2, ..., s_{(j-1)M+i-1}\}\) is not necessarily connected, i.e., $A^{\ell}_{(i,j)}$ cannot be defined by (\ref{defnbiggestbaseset}).
}

{\examp\label{eightnodesspace} Let \(S = \{i: 1\le i\le 8\}\) and \(s_i=i\) for \(1\le i\le 8\). The neighborhood of each \(s_i\) consist of 4 nodes which are connected to \(s_i\) by a dashed edge as illustrated by the following graph.

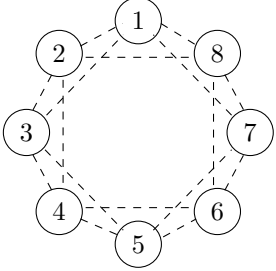
\begin{tikzpicture}
%\draw[help lines] (0,0) grid (2,2);
\begin{scope}[node distance=10mm]
\node[shape=circle,] (a) at (1,1) {};
\node[shape=circle,draw] [left=of a] {3};
\node[shape=circle,draw] [right=of a] {7};
\node[shape=circle,draw] [above=of a] {1};
\node[shape=circle,draw] [below=of a] {5};
\node[shape=circle,draw] [above left=of a] {2};
\node[shape=circle,draw] [above right=of a] {8};
\node[shape=circle,draw] [below left=of a] {4};
\node[shape=circle,draw] [below right=of a] {6};
\draw [dashed] (0.23,2.15) -- (0.79,2.45); % 1 - 2
\draw [dashed] (-0.4,1.3) -- (-0.15,1.77); % 2 - 3
\draw [dashed] (-0.4,0.7) -- (-0.15,0.23); % 3 - 4
\draw [dashed] (0.23,-0.15) -- (0.79,-.39); % 4 - 5
\draw [dashed] (1.33,-0.39) -- (1.77,-0.15); % 5 - 6
\draw [dashed] (2.15,0.23) -- (2.4,0.7); % 6 - 7
\draw [dashed] (2.15,1.77) -- (2.4,1.3); % 7 - 8
\draw [dashed] (1.3,2.45) -- (1.79,2.15); % 8 - 1

\draw [dashed] (0.0,1.7) -- (0.0,0.3); % 2 - 4
\draw [dashed] (0.3,0.0) -- (1.7,0.0); % 4 - 6
\draw [dashed] (2.0,0.3) -- (2.0,1.7); % 6 - 8
\draw [dashed] (0.3,2) -- (1.7,2); % 8 - 2

\draw [dashed] (-0.22,1.23) -- (0.80,2.30); % 1 - 3
\draw [dashed] (-0.24,0.77) -- (0.80,-0.30); % 3 - 5
\draw [dashed] (1.2,-0.3) -- (2.24,0.77); % 5 - 7
\draw [dashed] (2.24,1.23) -- (1.2,2.3); % 7 - 1

\end{scope}
\end{tikzpicture}

It can be verified that the neighborhood \(\partial(s_i)\) of each \(s_i\) as a subset of \(S\) is connected. For example, \(\partial(s_1)=\partial(1)=\{2,3,7,8\}\) is connected. But note that \(\partial(s_7)\bigcap\{s_i, 1\le i\le 6\}=\{1,5,6\}\) is not connected. So we cannot define base set \(A_{s_7}=\{1,5,6\}\).
}

{\rem Example \ref{eightnodesspace} shows that the condition that the neighborhood of each site is connected cannot imply \(\partial(s_i)\bigcap\{s_j, 1\le j\le i-1\}\) is connected for all \(2\le i\le N\). In fact, for abstract space \((S, \partial)\), it is difficult to find simple condition to guarantee that \(A_{s_i}=\partial(s_i)\bigcap\{s_j, 1\le j\le i-1\}\), \(\forall 2\le i\le N\). But fortunate enough, Example \ref{twodimensionalimagespace} illustrated that for many real random field applications, there exist natural neighborhood systems which have good properties and can satisfy our purpose very well.
}

\section{Novel Algorithm for Simulating Random Fields} \label{se:Algor_Sim_Ran_Field}

Let \((S,\partial)\) be a connected space of sites, with valid setup
\(\{(s_i, A_{s_i})\}_{i=1}^N\).
 For \(1\le i\le N\), \(d_i\in \mathbb{N}\) is the cardinality of \(\mathbf{X}_{s_i}\),
we denote \(\mathbf{X}_{s_i}=\{x_{s_i}^1, ..., x_{s_i}^{d_i}\}\).
Based on Proposition \ref{mainprop:alg}, we have a novel algorithm for
simulating $X$ with state space $\mathcal{X}$, given marginal
probabilities $\{\pi_{s_i}(x_{s_i}): x_{s_i}\in \mathbf{X}_{s_i}\}_{i=1}^N$
and given covariances of nearby sites \(\{ \beta_{s_i,t_i}: t_i\in A_{s_i} \}\), for \(2\le i\le N\).

Do for \(i=1,\dots, N\):
\begin{enumerate}
\item Base on Definition \ref{SequentiallySimulatable}, compute
\begin{eqnarray*}
   \Pi(X_{s_{i-1}}=x_{s_{i-1}}, ...,X_{s_1}=x_{s_1}) &=& \prod_{k=1}^{i-1} \Pi( X_{s_k}=x_{s_k}|X_{s_{k-1}}=x_{s_{k-1}}, ...,X_{s_{1}}=x_{s_{1}} ) \nonumber \\
   &=& \prod_{k=1}^{i-1} \Pi( X_{s_k}=x_{s_k}|X_{A_{s_k}}=x_{A_{s_k}})
    \nonumber
\end{eqnarray*}
for all chosen combinations of \(x_{s_1}\), ..., \(x_{s_{i-1}}\). Here we choose \(x_{s_1}\), ..., \(x_{s_{i-1}}\) as follows: for each \(1\le k\le i-1\), if \(s_k\in A_{s_i}\), we use the simulated \(x_{s_k}\); otherwise, we enumerate \(x_{s_k} \in \mathbf{X}_{s_k}\).

\item Take marginal to get \( \Pi(X_{A_{s_i}}=x_{A_{s_i}})\):
\begin{equation}
 \Pi(X_{A_{s_i}}=x_{A_{s_i}}) = \sum_{ s_k\notin A_{s_i}, 1\le k\le i-1} \Pi(X_{s_{i-1}}=x_{s_{i-1}}, ...,X_{s_1}=x_{s_1}). \nonumber
\end{equation}

\item Based on \(\Pi(X_{A_{s_i}}=x_{A_{s_i}})\), compute
\( \Pi(\left. X_{s_i}=x_{s_i}^{j}\right| X_{A_{s_i}}=x_{A_{s_i}}) \) for \( 1\le j\le d_i \), using (\ref{CondProb}).

\item Generate a \([0,1]\)-uniform random variable \(U\). For the given \(U\), there exists unique \( 1\le j\le d_i \) such that
\[
\sum_{u=1}^{j-1} \Pi(\left. X_{s_i}=x_{s_i}^{u}\right| X_{A_{s_i}}=x_{A_{s_i}}) \leq U < \sum_{u=1}^{j} \Pi(\left. X_{s_i}=x_{s_i}^{u}\right| X_{A_{s_i}}=x_{A_{s_i}}).
\]
Then set \(X_{s_i}= x_{s_i}^{j}\). For notational convenience, we supress superscript \(j\)\ and use \(x_{s_i}\) to indicate the simulated value \(x_{s_i}^{j}\) of \(X_{s_i}\).
\end{enumerate}

{\rem\label{Re:Algor_Sim_Ran_Field} There exists an alternative, more efficient way of computing \(\Pi(X_{A_{s_i}}=x_{A_{s_i}})\): View \(A_{s_i}\) as \(S\), and let \(\{(t_j,B_{t_j})\}_{j=1}^{n_i}\) be a valid setup for \(A_{s_i}\), where \(n_i\) is the number of sites in \(A_{s_i}\). Now, compute \(\Pi(X_{A_{s_i}}=x_{A_{s_i}})\) as follows:
\[
  \Pi(X_{A_{s_i}}=x_{A_{s_i}}) = \prod_{j=1}^{n_i} \Pi( X_{t_j}=x_{t_j}|X_{B_{t_j}}=x_{B_{t_j}}).
\]
In Example \ref{example1}, we illustrated how to compute \(\Pi(X_{A_5}=x_{A_5})\) by treating \(A_5\) as \(S\). }

\section{Summary of Applications of Novel Algorithm }\label{se:Appl_Novel_Algorithm_to_Image_Analysis}

In this section, we summarize our prior applications of simulation
algorithm of Section \ref{se:Algor_Sim_Ran_Field} to image generation, data authentication and target recognition.

\subsection{Application to Image Generation}
\label{Image_Generation}

Generating KNW-CAPTCHAs is an application of our simulation algorithm to image analysis. CAPTCHA is the acronym for ``Completely Automated Public Turing test to tell Computers and Hummans Apart'' (see \cite{blum-captcha}) and is widely used to prevent online resources intended for humans from abuse by automated agents. CAPTCHAs often appear to be images of characters or digits designed easy to read by humans and difficult to crack by computer programs. Kouritzin, Newton and Wu \cite{KoNeWu} proposed a novel method for generating a type of CAPTCHAs named ``KNW-CAPTCHAs'' through random field simulation stated in Section \ref{se:Algor_Sim_Ran_Field} with common pixel state space \(\mathbf{X}=\{1,-1\}=\{black, white\}\). Roughly speaking, KNW-CAPTCHAs are generated by specifying proper pixel marginal probabilities and pixel-pixel covariances of alphabets and embedding these quantities into KNW conditional probabilities - formula (\ref{Captchaformula}). One predominance of this method is that these random field CAPTCHAs can be simulated in real time, yet another is significant resistance of these CAPTCHAs to attack. More details on this application can be found in \cite{KoNeWu}. With the algorithm in Section \ref{se:Algor_Sim_Ran_Field} and the alternative algorithm for generating Markov random field at the end of Section \ref{se:AlgrthmConstrnts}, we are also able to generate many variants of KNW-CAPTCHAs.

\subsection{Data Authentication Application}\label{Data_Authentication}

Data authentication is classifying data as true or fabricated and has been used in fraud detection and verification of data samples including financial data \cite{HILLTP}. Detecting fake coin flip sequence was an application of our simulation algorithm to data authentication, see Kouritzin et al  \cite{Kouritzin2008} where site space \(S=\{1,2, ..., N\}\) consists of time units and coin state space is \(\mathbf{X}=\{1,-1\}=\{head, tail\}\). In \cite{Kouritzin2008}, a filtering method was applied to simplified fraud-detection problem, classifying coin flip sequences as either ``faked" (i.e., generated by a person) or ``real" (i.e., generated by perfect flipping of a true coin). A true coin flip sequence has the expected head counts approximating \(\frac{1}{2}\) of total flips. For this reason, a faker that deviates from expected behaviour in one time period has to be compensated by later deviating from the expected behaviour in an opposite way such that deviant behaviour averages out. Consequently, marginal probabilities and pair-wise covariances between each flip and the flips that preceded in time were used to describe faker's behaviours. It followed that simulating fake coin flip sequences in real time became a fundamental step of solving the simplified fraud-detection problem. The much simpler but preceded version of the algorithm in Section \ref{se:Algor_Sim_Ran_Field} was developed to simulate the faked coin flip sequences which are just specific random fields on the set \(S\) with state pace \(\mathbf{X}=\{head, tail\}\).

\subsection{Target Recognition Application}\label{Target_Recognition}
Analogue of simulation algorithm in Section \ref{se:Algor_Sim_Ran_Field} was also applied to target recognition by Kouritzin, Luo, Newton and Wu \cite{Kouritzin2009}. Kouritzin et al considered an imaginary detection problem of hidden targets such as rocket launchers in random forest. The rocket lauchers sit still in forest, and a surveillance aircraft or unmanned vehicle equipped with an electro-optic camera is flying over the forest to capture the forest with hiden target rocket launchers. The camera cannot penetrate through the foliage, and the foliage blocks the latent ground objects from detection. However, the camera can observe partial images of the ground through the gaps among the leaves. The forest is random in the sense that the foliage coverage is encoded with a binary representation (i.e., foliage and no foliage), and is correlated in adjacent regions. Similarly, the ground is a mixture of the grass and soil, and its color is randomly either green or brown with some type of correlation structure. The rocket launchers are camouflaged with the colors of grass (green) and soil (brown) and have another type of correlation structure. The information from within the forest called the observations, can be easily obtained from stored historical overhead pictures and analyzed pixel by pixel for each small area. As the captcha generation, simulation algorithm similar to that in Section \ref{se:Algor_Sim_Ran_Field} was used to generate the image of random forest, ground and hidden targets. The weight functions used by the SERP (selectively resampling particle) filter \cite{Ballantyne2002} resort to Proposition 2.7 of \cite{Kouritzin2009}, which is a specific vector version of Proposition \ref{mainprop:alg}.

\section{Algorithm Properties and Constraints}\label{se:AlgrthmConstrnts}

An effective algorithm that can generate random fields is related to two properties: regularity and consistency. Regularity means that the right hand side of (\ref{CondProb}) is within [0,1]. Consistency means the permutation and marginality properties. In this section, we explore the constraints that can guarantee these properties. The outline of this section is as follows. First, we give the necessary and sufficient conditions of the regularity type for singleton case of \(A_{s_i}\) in Proposition \ref{prop:condition}. Second, we consider conditions for permutation and marginality properties for \(A\subseteq S\) where sites in \(A\) are neighbors of each other. For random field \(X_A\), there are no conditions needed for the marginality property, and there exists a necessary and sufficient condition presented in Proposition \ref{permutation:necesuffi} and \ref{permutation:necesuffm} for the permutation property. Third, for an arbitrary site space \((S,\partial)\), the neighborhood system \(\partial\) is extended to \(\partial'\) such that each pair of sites are neighbors of each other under \(\partial'\), i.e., for each site \(s\in S\), let \(\partial'(s) = S\setminus\{s\}\). Furthermore, the covariances between each site \(s\in S\) and sites outside of \(\partial(s)\) are assumed to be zero. Based on these assumptions, we propose an alternative algorithm that can generate a true Markov random field with the remarkable property of consistency.

We first determine the conditions of regularity type, i.e., the constraints for \(\{ \beta_{s_i,t_i}: t_i\in A_{s_i} \}_{i=1}^N\) of Proposition \ref{mainprop:alg} that cause it to produce a value in \([0,1]\). A theoretical constraint relating a particular \(\beta_{s_i,t_i}\) to the given marginal probabilities $\{\pi_{s_i}(x_{s_i}): x_{s_i}\in \mathbf{X}_{s_i}\}$ and $\{\pi_{t_i}(x_{t_i}): x_{t_i}\in \mathbf{X}_{t_i}\}$ is given in the following proposition.
\begin{prop}
\label{prop:condition}
If \(S=\{s_i\}_{i=1}^N\) and \( A_{s}=\{t\} \) is a singleton, then a necessary and sufficient condition on \(\beta_{s,t}\) for (\ref{CondProb}) to be in [0,1] is that  \(\beta_{s,t}\) is a valid covariance; which corresponds to the condition
\begin{eqnarray}\label{CondEqn}
 \beta_{s,t}&\in & \biggl[\max_{\substack{x_{s} \in \mathbf{X}_{s}, x_{s}\not=\mu_{\tilde{\pi}_{s}}\\ x_{t}\in \mathbf{X}_{t}, x_{t}\not=\mu_{\tilde{\pi}_{t}}}}\biggl(\frac{-\sigma^2_{\tilde{\pi}_{s}}\sigma^2_{\tilde{\pi}_{t}}\pi_{s}(x_{s}) \pi_{t}(x_{t})}
{\tilde{\pi}_{s}(x_{s}) \tilde{\pi}_{t}(x_{t})(x_{s}-\mu_{\tilde{\pi}_{s}})(x_{t}-\mu_{\tilde{\pi}_{t}})}\bigwedge \frac{\sigma^2_{\tilde{\pi}_{s}}\sigma^2_{\tilde{\pi}_{t}}(1-\pi_{s}(x_{s})) \pi_{t}(x_{t})}
{\tilde{\pi}_{s}(x_{s}) \tilde{\pi}_{t}(x_{t})(x_{s}-\mu_{\tilde{\pi}_{s}})(x_{t}-\mu_{\tilde{\pi}_{t}})}\biggl), \nonumber \\
&& \min_{\substack{x_{s} \in \mathbf{X}_{s}, x_{s}\not=\mu_{\tilde{\pi}_{s}}\\ x_{t}\in \mathbf{X}_{t}, x_{t}\not=\mu_{\tilde{\pi}_{t}}}}\biggl(\frac{-\sigma^2_{\tilde{\pi}_{s}}\sigma^2_{\tilde{\pi}_{t}}\pi_{s}(x_{s}) \pi_{t}(x_{t})}
{\tilde{\pi}_{s}(x_{s}) \tilde{\pi}_{t}(x_{t})(x_{s}-\mu_{\tilde{\pi}_{s}})(x_{t}-\mu_{\tilde{\pi}_{t}})}\bigvee \frac{\sigma^2_{\tilde{\pi}_{s}}\sigma^2_{\tilde{\pi}_{t}}(1-\pi_{s}(x_{s})) \pi_{t}(x_{t})}
{\tilde{\pi}_{s}(x_{s}) \tilde{\pi}_{t}(x_{t})(x_{s}-\mu_{\tilde{\pi}_{s}})(x_{t}-\mu_{\tilde{\pi}_{t}})}\biggl)
\biggl]. \nonumber \\
\end{eqnarray}
\end{prop}

\proof For \(x_{s}\in \mathbf{X}_{s}, x_{t}\in \mathbf{X}_{t} \), we have by (\ref{CondProb}) that
\begin{eqnarray}\label{CondiEquiv}
 &&\Pi(X_{s}=x_{s}| X_{t}=x_{t})\in [0,1] \nonumber \\
 &\Leftrightarrow & \frac{\sigma^2_{\tilde{\pi}_{s}}\sigma^2_{\tilde{\pi}_{t}}\pi_{s}(x_{s}) \pi_{t}(x_{t})+\tilde{\pi}_{s}(x_{s}) \tilde{\pi}_{t}(x_{t})(x_{s}-\mu_{\tilde{\pi}_{s}})\beta_{s,t}(x_{t}-\mu_{\tilde{\pi}_{t}})}{\sigma^2_{\tilde{\pi}_{s}}\sigma^2_{\tilde{\pi}_{t}} \pi_{t}(x_{t})}\in [0,1] \nonumber \\
 &\Leftrightarrow & 0\le \sigma^2_{\tilde{\pi}_{s}}\sigma^2_{\tilde{\pi}_{t}}\pi_{s}(x_{s}) \pi_{t}(x_{t})+\tilde{\pi}_{s}(x_{s}) \tilde{\pi}_{t}(x_{t})(x_{s}-\mu_{\tilde{\pi}_{s}})\beta_{s,t}(x_{t}-\mu_{\tilde{\pi}_{t}})\le \sigma^2_{\tilde{\pi}_{s}}\sigma^2_{\tilde{\pi}_{t}} \pi_{t}(x_{t})\nonumber \\
  &\Leftrightarrow & -\sigma^2_{\tilde{\pi}_{s}}\sigma^2_{\tilde{\pi}_{t}}\pi_{s}(x_{s}) \pi_{t}(x_{t}) \le \tilde{\pi}_{s}(x_{s}) \tilde{\pi}_{t}(x_{t})(x_{s}-\mu_{\tilde{\pi}_{s}})\beta_{s,t}(x_{t}-\mu_{\tilde{\pi}_{t}})\le \sigma^2_{\tilde{\pi}_{s}}\sigma^2_{\tilde{\pi}_{t}} (1-\pi_{s}(x_{s})) \pi_{t}(x_{t})\nonumber \\
\end{eqnarray}
By (\ref{CondiEquiv}), there is no constraint on \(\beta_{s,t}\) when \(x_{s}=\mu_{\tilde{\pi}_{s}}\) or \(x_{t}=\mu_{\tilde{\pi}_{t}}\). Therefore,  \(\beta_{s,t}\) being consistent with a proper right hand side of (\ref{CondProb}) is equivalent to,
\begin{eqnarray}\label{CondiEquiv2}
 \beta_{s,t} &\in& \biggl[\frac{-\sigma^2_{\tilde{\pi}_{s}}\sigma^2_{\tilde{\pi}_{t}}\pi_{s}(x_{s}) \pi_{t}(x_{t})}
{\tilde{\pi}_{s}(x_{s}) \tilde{\pi}_{t}(x_{t})(x_{s}-\mu_{\tilde{\pi}_{s}})(x_{t}-\mu_{\tilde{\pi}_{t}})}\bigwedge \frac{\sigma^2_{\tilde{\pi}_{s}}\sigma^2_{\tilde{\pi}_{t}}(1-\pi_{s}(x_{s})) \pi_{t}(x_{t})}
{\tilde{\pi}_{s}(x_{s}) \tilde{\pi}_{t}(x_{t})(x_{s}-\mu_{\tilde{\pi}_{s}})(x_{t}-\mu_{\tilde{\pi}_{t}})}, \nonumber \\
&& \frac{-\sigma^2_{\tilde{\pi}_{s}}\sigma^2_{\tilde{\pi}_{t}}\pi_{s}(x_{s}) \pi_{t}(x_{t})}
{\tilde{\pi}_{s}(x_{s}) \tilde{\pi}_{t}(x_{t})(x_{s}-\mu_{\tilde{\pi}_{s}})(x_{t}-\mu_{\tilde{\pi}_{t}})} \bigvee \frac{\sigma^2_{\tilde{\pi}_{s}}\sigma^2_{\tilde{\pi}_{t}}(1-\pi_{s}(x_{s})) \pi_{t}(x_{t})}
{\tilde{\pi}_{s}(x_{s}) \tilde{\pi}_{t}(x_{t})(x_{s}-\mu_{\tilde{\pi}_{s}})(x_{t}-\mu_{\tilde{\pi}_{t}})}
\biggl]
\end{eqnarray}
for all \(x_{s} \in \mathbf{X}_{s}, x_{s}\not=\mu_{\tilde{\pi}_{s}}\), and \(x_{t}\in \mathbf{X}_{t}, x_{t}\not=\mu_{\tilde{\pi}_{t}}\). From (\ref{CondiEquiv2}), (\ref{CondEqn}) follows.

In the Appendix, we will prove that
\begin{equation}\label{Jointprob}
\Pi(X_{s}=x_{s},  X_{t}=x_{t}) = \pi_{s}(x_{s})\pi_{t}(x_{t}) + \frac{\displaystyle{\tilde{\pi}_{s}(x_{s})(x_{s}-\mu_{\tilde{\pi}_{s}})cov(X_{s},X_{t})\tilde{\pi}_{t}(x_{t})(x_{t}-\mu_{\tilde{\pi}_{t}})}}
{\displaystyle{\sigma^2_{\tilde{\pi}_{s}}\displaystyle{\sigma^2_{\tilde{\pi}_{t}}}}}
\end{equation}
By multiplying both sides of (\ref{Jointprob}) by \(\frac{\sigma^2_{\tilde{\pi}_{s}}\sigma^2_{\tilde{\pi}_{t}}}{\tilde{\pi}_{s}(x_{s})\tilde{\pi}_{t}(x_{t})}\), for any \(x_{s} \in \mathbf{X}_{s}, x_{t}\in \mathbf{X}_{t}\), we can verify that
\begin{eqnarray}\label{betaValidCov}
  (x_{s}-\mu_{\tilde{\pi}_{s}})cov(X_{s},X_{t})(x_{t}-\mu_{\tilde{\pi}_{t}}) & = &\frac{\sigma^2_{\tilde{\pi}_{s}}\sigma^2_{\tilde{\pi}_{t}}}{\tilde{\pi}_{s}(x_{s})\tilde{\pi}_{t}(x_{t})}[\Pi(X_{s}=x_{s}, X_{t}=x_{t})\Pi(X_{s}\not=x_{s}, X_{t}\not=x_{t}) \nonumber \\
    &{}& -\Pi(X_{s}=x_{s}, X_{t}\not=x_{t})\Pi(X_{s}\not=x_{s}, X_{t}=x_{t})].
\end{eqnarray}

It follows that
\begin{eqnarray}\label{ConstraintIneq}
-\frac{\sigma^2_{\tilde{\pi}_{s}}\sigma^2_{\tilde{\pi}_{t}}\pi_{s}(x_{s})\pi_{t}(x_{t})}{\tilde{\pi}_{s}(x_{s})\tilde{\pi}_{t}(x_{t})}  &\le & -\frac{\sigma^2_{\tilde{\pi}_{s}}\sigma^2_{\tilde{\pi}_{t}}}{\tilde{\pi}_{s}(x_{s})\tilde{\pi}_{t}(x_{t})} \Pi(X_{s}=x_{s}, X_{t}\not=x_{t})\Pi(X_{s}\not=x_{s}, X_{t}=x_{t}) \nonumber \\
   &\le & (x_{s}-\mu_{\tilde{\pi}_{s}})cov(X_{s},X_{t}) (x_{t}-\mu_{\tilde{\pi}_{t}})  \\
   &\le & \frac{\sigma^2_{\tilde{\pi}_{s}}\sigma^2_{\tilde{\pi}_{t}}}{\tilde{\pi}_{s}(x_{s})\tilde{\pi}_{t}(x_{t})} \Pi(X_{s}=x_{s}, X_{t}=x_{t})\Pi(X_{s}\not=x_{s}, X_{t}\not=x_{t}) \nonumber \\
   &\le & \frac{\sigma^2_{\tilde{\pi}_{s}}\sigma^2_{\tilde{\pi}_{t}}(1-\pi_{s}(x_{s}))\pi_{t}(x_{t})}{\tilde{\pi}_{s}(x_{s})\tilde{\pi}_{t}(x_{t})}. \nonumber
\end{eqnarray}

The inequalities in (\ref{ConstraintIneq}) are actually tight. If for \(x_{s} \in \mathbf{X}_{s}\) and \(x_{t}\in \mathbf{X}_{t}\), \(\{X_{s}=x_{s}\}= \{X_{t}=x_{t}\}\), then \(\{X_{s}\not=x_{s}\} = \{X_{t}\not=x_{t}\}\). It follows that
\[
  \{X_{s}=x_{s}, X_{t}=x_{t}\} = \{X_{t}=x_{t}\}, \hspace{4mm}
  \{X_{s}\not=x_{s}, X_{t}\not=x_{t}\} = \{X_{s}\not=x_{s}\},
\]
and
\[
  \{X_{s}=x_{s}, X_{t}\not=x_{t}\} =\emptyset.
\]
Then, by (\ref{betaValidCov}), we have
\begin{eqnarray*}
   (x_{s}-\mu_{\tilde{\pi}_{s}})cov(X_{s},X_{t}) (x_{t}-\mu_{\tilde{\pi}_{t}})& = & \frac{\sigma^2_{\tilde{\pi}_{s}}\sigma^2_{\tilde{\pi}_{t}}}{\tilde{\pi}_{s}(x_{s})\tilde{\pi}_{t}(x_{t})} [\Pi(X_{t}=x_{t})\Pi(X_{s}\not=x_{s})-0\times\Pi(X_{s}\not=x_{s}, X_{t}=x_{t})] \\
   &=& \frac{\sigma^2_{\tilde{\pi}_{s}}\sigma^2_{\tilde{\pi}_{t}}(1-\pi_{s}(x_{s}))\pi_{t}(x_{t})}{\tilde{\pi}_{s}(x_{s})\tilde{\pi}_{t}(x_{t})}.
\end{eqnarray*}
Similarly, if \(\{X_{t}\not=x_{t}\} = \{X_{s}=x_{s}\}\), we have
\[
   -\frac{\sigma^2_{\tilde{\pi}_{s}}\sigma^2_{\tilde{\pi}_{t}}\pi_{s}(x_{s})\pi_{t}(x_{t})}{\tilde{\pi}_{s}(x_{s})\tilde{\pi}_{t}(x_{t})} = (x_{s}-\mu_{\tilde{\pi}_{s}})cov(X_{s},X_{t})(x_{t}-\mu_{\tilde{\pi}_{t}}).
\]
Therefore, by (\ref{CondiEquiv}) and (\ref{ConstraintIneq}), (\ref{CondEqn}) is exactly the constraint required for \(\beta_{s,t}\) to be a valid covariance.
\endproof

Secondly, we consider the case when \( A_{s} \) is not a singleton. For this case, we assume that \(cov(X_{s}, X_{t})=\beta_{s}\) for every \(t\in A_{s} \), i.e., the state \(X_{s}\) on the site \(s\) has same covariance with each state in \( A_{s} \). A simple necessary condition on \(\beta_{s}\) is given in the following corollary.

\begin{cor}
\label{corllary:condition2}
Assume that \(cov(X_{s}, X_{t})=\beta_{s}\) for each \(t \in A_s\). Then, a necessary condition on \(\beta_{s}\) for (\ref{CondProb}) to be in [0,1] is
\begin{eqnarray}\label{CondEqn2}
\beta_s &\in &\biggl[\max_{\substack{x_{s}\in \mathbf{X}_{s}, x_{s}\not=\mu_{\tilde{\pi}_{s}},\\ x_t\in \mathbf{X}_t, x_t\not=\mu_{\tilde{\pi}_t}, \\ t\in A_s}}\biggl(\frac{-\sigma^2_{\tilde{\pi}_{s}}\sigma^2_{\tilde{\pi}_t} \pi_{s}(x_{s}) \pi_t(x_t)}
{{\tilde{\pi}_{s}(x_{s})\tilde{\pi}_t(x_t)}(x_{s}-\mu_{\tilde{\pi}_{s}})(x_t-\mu_{\tilde{\pi}_t})}\bigwedge \frac{\sigma^2_{\tilde{\pi}_{s}}\sigma^2_{\tilde{\pi}_t} (1-\pi_{s}(x_{s})) \pi_t(x_t)}
{{\tilde{\pi}_{s}(x_{s})\tilde{\pi}_t(x_t)}(x_{s}-\mu_{\tilde{\pi}_{s}})(x_t-\mu_{\tilde{\pi}_t})}\biggl),\nonumber \\
&& \min_{\substack{x_{s}\in \mathbf{X}_{s}, x_{s}\not=\mu_{\tilde{\pi}_{s}},\\ x_t\in \mathbf{X}_t, x_t\not=\mu_{\tilde{\pi}_t}, \\ t\in A_s}}\biggl(\frac{-\sigma^2_{\tilde{\pi}_{s}}\sigma^2_{\tilde{\pi}_t} \pi_{s}(x_{s}) \pi_t(x_t)}
{{\tilde{\pi}_{s}(x_{s})\tilde{\pi}_t(x_t)}(x_{s}-\mu_{\tilde{\pi}_{s}})(x_t-\mu_{\tilde{\pi}_t})}\bigvee \frac{\sigma^2_{\tilde{\pi}_{s}}\sigma^2_{\tilde{\pi}_t} (1-\pi_{s}(x_{s})) \pi_t(x_t)}
{{\tilde{\pi}_{s}(x_{s})\tilde{\pi}_t(x_t)}(x_{s}-\mu_{\tilde{\pi}_{s}})(x_t-\mu_{\tilde{\pi}_t})}\biggl)
\biggl].\nonumber \\
\end{eqnarray}
\end{cor}

\proof For each \(t\in A_s\), we get a condition on \(\beta_{s}\) by applying Proposition \ref{prop:condition} to \(A_{s,t}=\{t\}\). (\ref{CondEqn2}) follows by combining all these conditions into one. This is also the reason that (\ref{CondEqn2}) is only necessary, not sufficient for (\ref{CondProb}) to be in [0,1].   \endproof

Assume that \(A=\{s_i: 1\le i\le n\}\subseteq S\) where sites in \(A\) are neighbors of each other under \(\partial\), we denote the joint probability mass function of \(X_{s_{1}}\), ...,\(X_{s_{n}}\) by \(\Pi_{s_1, ..., s_n}(x_{s_{1}}, ..., x_{s_{n}})\) which by formula (\ref{CondProb}) satisfies:
\begin{eqnarray}\label{recursivejointpmf}
 && \Pi_{s_1, ..., s_i}(x_{s_{1}}, ..., x_{s_{i}})  \\
 &=& \Pi_{s_1, ..., s_{i-1}}(x_{s_{1}}, ..., x_{s_{i-1}}) \pi_{s_i}(x_{s_i}) + \frac{\tilde{\pi}_{s_i}(x_{s_i})(x_{s_i}-\mu_{\tilde{\pi}_{s_i}})}{\sigma^2_{\tilde{\pi}_{s_i}}}\sum_{j=1}^{i-1} \prod_{k=1,k\not=j}^{i-1}\hat{\pi}_{s_k}(x_{s_k}) \cdot \frac{\beta_{s_i,s_j}\tilde{\pi}_{s_j}(x_{s_j})(x_{s_j}-\mu_{\tilde{\pi}_{s_j}})}{\sigma^2_{\tilde{\pi}_{s_j}}} \nonumber
\end{eqnarray}
for \(i=n, n-1, ..., 3\) and
\begin{equation}\label{jointpmfs1s2}
 \Pi_{s_1, s_2}(x_{s_{1}}, x_{s_{2}}) = \pi_{s_2}(x_{s_2}) \pi_{s_1}(x_{s_1}) + \frac{\tilde{\pi}_{s_2}(x_{s_2})(x_{s_2}-\mu_{\tilde{\pi}_{s_2}})\beta_{s_2,s_1}\tilde{\pi}_{s_1}(x_{s_1})(x_{s_1}-\mu_{\tilde{\pi}_{s_1}})}{\sigma^2_{\tilde{\pi}_{s_2}}\sigma^2_{\tilde{\pi}_{s_1}}}.
\end{equation}
We have the following lemma giving the closed form of \(\Pi_{s_1, ..., s_n}(x_{s_{1}}, ..., x_{s_{n}})\):
\begin{lem}\label{lem:jointpmf} Based on the assumptions on \(A\), the joint probability mass function of \(X_{s_1}, ..., X_{s_n}\) is
\begin{eqnarray}\label{closedform_neighbors}
 && \Pi_{s_1, ..., s_n}(x_{s_{1}}, ..., x_{s_{n}}) \nonumber \\
 &=& \sum_{i=2}^n \biggl[\biggl(\frac{\tilde{\pi}_{s_i}(x_{s_i})(x_{s_i}-\mu_{\tilde{\pi}_{s_i}})}{\sigma^2_{\tilde{\pi}_{s_i}}}\sum_{j=1}^{i-1} \prod_{k=1,k\not=j}^{i-1}\hat{\pi}_{s_k}(x_{s_k}) \cdot \frac{\beta_{s_i,s_j}\tilde{\pi}_{s_j}(x_{s_j})(x_{s_j}-\mu_{\tilde{\pi}_{s_j}})}{\sigma^2_{\tilde{\pi}_{s_j}}}\biggl)\cdot \prod_{k=i+1}^n \pi_{s_k}(x_{s_k}) \biggl] + \prod_{i=1}^n \pi_{s_i}(x_{s_i}) \nonumber \\
\end{eqnarray}
for each \(x_{s_i} \in \mathbf{X}_{s_i}\) (\(1\le i\le n\)).
\end{lem}

\proof Repeated use of (\ref{recursivejointpmf}) followed by (\ref{jointpmfs1s2}) yields
\begin{eqnarray}
 && \Pi_{s_1, ..., s_n}(x_{s_{1}}, ..., x_{s_{n}})  \nonumber \\
 &=& \Pi_{s_1, ..., s_{n-1}}(x_{s_{1}}, ..., x_{s_{n-1}}) \pi_{s_n}(x_{s_n}) + \frac{\tilde{\pi}_{s_n}(x_{s_n})(x_{s_n}-\mu_{\tilde{\pi}_{s_n}})}{\sigma^2_{\tilde{\pi}_{s_n}}}\sum_{j=1}^{n-1} \prod_{k=1,k\not=j}^{n-1}\hat{\pi}_{s_k}(x_{s_k}) \cdot \frac{\beta_{s_n,s_j}\tilde{\pi}_{s_j}(x_{s_j})(x_{s_j}-\mu_{\tilde{\pi}_{s_j}})}{\sigma^2_{\tilde{\pi}_{s_j}}} \nonumber \\
 &=& \Pi_{s_1, ..., s_{n-2}}(x_{s_{1}}, ..., x_{s_{n-2}}) \pi_{s_{n-1}}(x_{s_{n-1}})\pi_{s_n}(x_{s_n}) \nonumber \\
 && + \frac{\tilde{\pi}_{s_{n-1}}(x_{s_{n-1}})(x_{s_{n-1}}-\mu_{\tilde{\pi}_{s_{n-1}}})}{\sigma^2_{\tilde{\pi}_{s_{n-1}}}}\sum_{j=1}^{n-2} \prod_{k=1,k\not=j}^{n-2}\hat{\pi}_{s_k}(x_{s_k}) \cdot \frac{\beta_{s_{n-1},s_j}\tilde{\pi}_{s_j}(x_{s_j})(x_{s_j}-\mu_{\tilde{\pi}_{s_j}})}{\sigma^2_{\tilde{\pi}_{s_j}}} \cdot \pi_{s_n}(x_{s_n})  \nonumber \\
 && + \frac{\tilde{\pi}_{s_n}(x_{s_n})(x_{s_n}-\mu_{\tilde{\pi}_{s_n}})}{\sigma^2_{\tilde{\pi}_{s_n}}}\sum_{j=1}^{n-1} \prod_{k=1,k\not=j}^{n-1}\hat{\pi}_{s_k}(x_{s_k}) \cdot \frac{\beta_{s_n,s_j}\tilde{\pi}_{s_j}(x_{s_j})(x_{s_j}-\mu_{\tilde{\pi}_{s_j}})}{\sigma^2_{\tilde{\pi}_{s_j}}} \nonumber \\
&=& \sum_{i=2}^n \biggl[\biggl(\frac{\tilde{\pi}_{s_i}(x_{s_i})(x_{s_i}-\mu_{\tilde{\pi}_{s_i}})}{\sigma^2_{\tilde{\pi}_{s_i}}}\sum_{j=1}^{i-1} \prod_{k=1,k\not=j}^{i-1}\hat{\pi}_{s_k}(x_{s_k}) \cdot \frac{\beta_{s_i,s_j}\tilde{\pi}_{s_j}(x_{s_j})(x_{s_j}-\mu_{\tilde{\pi}_{s_j}})}{\sigma^2_{\tilde{\pi}_{s_j}}}\biggl)\cdot \prod_{k=i+1}^n \pi_{s_k}(x_{s_k}) \biggl] + \prod_{i=1}^n \pi_{s_i}(x_{s_i}). \nonumber \\
\end{eqnarray}
\qed

{\rem\label{rem:lem:jointpmf} Our KNW quick simulation random field \(X_A\) on \(A\) satisfies the \textbf{marginality} condition defined in Definition \ref{ConsistencyCondi}. This can be proved by formula (\ref{closedform_neighbors}) and the facts that \(\displaystyle\sum_{x_{s_i} \in \mathbf{X}_{s_i}}\tilde{\pi}_{s_i}(x_{s_i})(x_{s_i}-\mu_{\tilde{\pi}_{s_i}})=0\), \(\displaystyle\sum_{x_{s_i} \in \mathbf{X}_{s_i}}\hat{\pi}_{s_i}(x_{s_i})=1\), \(\displaystyle\sum_{x_{s_i} \in \mathbf{X}_{s_i}}\pi_{s_i}(x_{s_i})=1\) for \(1\le i\le n\). Therefore, there are no derived conditions on \(\beta_{s_i, s_j}\) (\(1\le i\not=j\le n\)) and \(\pi_{s_i}\) (\(1\le i\le n\)) for \(X_A\) to satisfy the marginality condition.
}

Next, we determine conditions of permutation type. Assume that \(s,t\in S\) are neighbors of each other. We have by (\ref{CondProb}) and the multiplication rule that
\begin{equation}\label{consistency2Variables}
  \Pi(X_s=x_s|X_t=x_t)\Pi(X_t=x_t) = \Pi(X_t=x_t|X_s=x_s)\Pi(X_s=x_s),
\end{equation}
for \(x_s \in \mathbf{X}_s\), \(x_t\in \mathbf{X}_t\). Hence, \(\Pi_{t, s}(x_{t}, x_{s}) = \Pi_{s, t}(x_{s}, x_{t})\), i.e., when we compute the joint probability of two simulated random variables \(X_s\), \(X_t\), the order of simulating \(X_s\) and \(X_t\) does not matter. It means that the specific form of (\ref{CondProb}) guarantees the permutation property of any two random variables with their sites being neighbors. This is a very desirable property of formula (\ref{CondProb}) when the site space \(S\) is one-dimensional, e.g., \(S\subset\mathbb{N}\).

When there are \(n\ge 3\) sites which are neighbors of each other, we have a necessary and sufficient condition for permutation property given by Propposition \ref{permutation:necesuffi} below. First, we introduce some notations and preliminary results. For each \(1\le i \le n\) and \(1\le u_i\le d_i\), let
\begin{equation}\label{defntildez}
 \tilde{z}_{s_i}^{u_i} = \frac{\tilde{\pi}_{s_i}(x_{s_i}^{u_i})(x_{s_i}^{u_i}-\mu_{\tilde{\pi}_{s_i}})}{\sigma_{\tilde{\pi}_{s_i}}},
\end{equation}
where \(x_{s_i}^{u_i} \in \mathbf{X}_{s_i} = \{x_{s_i}^1, ..., x_{s_i}^{d_i}\}\). In the proof of Propposition \ref{permutation:necesuffi}, if we use \(x_{s_i}\), omitting \(u_i\) for \(x_{s_i}^{u_i}\), then we do the same thing for \(\tilde{z}_{s_i}^{u_i}\), i.e., we use \(\tilde{z}_{s_i}\). For any duple \(1\le i\not=j\le n\), let
\begin{equation}\label{defntilderho}
 \tilde{\rho}_{s_i,s_j} = \frac{\beta_{s_i,s_j}}{\sigma_{\tilde{\pi}_{s_i}}\sigma_{\tilde{\pi}_{s_j}}}
 =\frac{\sigma_{\pi_{s_i}}\sigma_{\pi_{s_j}}}{\sigma_{\tilde{\pi}_{s_i}}\sigma_{\tilde{\pi}_{s_j}}}\rho_{s_i,s_j}.
\end{equation}

Note that \(\tilde{z}_{s_i}^{u_i}\) is the weighted standardization of the element \(x_{s_i}^{u_i}\) by the weight \(\tilde{\pi}_{s_i}(x_{s_i}^{u_i})\) and \(\tilde{\rho}_{s_i,s_j}\) is an adjusted correlation coefficient of \(X_{s_i}\) and \(X_{s_j}\). Without loss of generality, we assume that \(\rho_{s_1,s_2}\not=0\) and \(\rho_{s_1,s_n}\not=0\), i.e.,  \(\tilde{\rho}_{s_1,s_2}\not=0\) and \(\tilde{\rho}_{s_1,s_n}\not=0\). We introduce \(d\) (\(= \sum_{i=1}^n d_i\)) variables \(y_{s_i}^{u_i}\) and denote \(d\) components
\begin{equation}\label{blockwiseCanonicalformsolutionform1}
  \hat{y}_{s_i}^{u_i} = \hat{\pi}_{s_i}(x_{s_i}^{u_i}) - \pi_{s_i}(x_{s_i}^{u_i})
\end{equation}
for \(1\le u_i\le d_i, 1\le i\le n\).

For each triple \(1\le i,j,k \le n\) (\(i,j\) and \(k\) are distinct integers) and \(1\le u_i\le d_i\),  \(1\le u_j\le d_j\) and  \(1\le u_k\le d_k\), we formulate equations in \(y\) as follows:
\begin{equation}\label{permutationcondiTriple}
  \tilde{z}_{s_i}^{u_i} \tilde{z}_{s_j}^{u_j}\tilde{\rho}_{s_i,s_j} y_{s_k}^{u_k} = \tilde{z}_{s_i}^{u_i} \tilde{z}_{s_k}^{u_k} \tilde{\rho}_{s_i,s_k} y_{s_j}^{u_j}  = \tilde{z}_{s_j}^{u_j}\tilde{z}_{s_k}^{u_k} \tilde{\rho}_{s_j,s_k} y_{s_i}^{u_i}
\end{equation}
(\ref{permutationcondiTriple}) is equivalent to the linear equations of \textbf{Blockwise Canonical form}:
\begin{eqnarray}\label{blockwiseCanonicalform}
  \tilde{z}_{s_j}^{u_j}\tilde{\rho}_{s_j,s_k} y_{s_i}^{u_i} &= &  \tilde{z}_{s_i}^{u_i} \tilde{\rho}_{s_i,s_k} y_{s_j}^{u_j}  \nonumber \\
  \tilde{z}_{s_k}^{u_k} \tilde{\rho}_{s_j,s_k} y_{s_i}^{u_i} &= & \tilde{z}_{s_i}^{u_i}  \tilde{\rho}_{s_i,s_j} y_{s_k}^{u_k}.
\end{eqnarray}
The coefficient matrix of the equations (\ref{blockwiseCanonicalform}) is
\begin{equation}\label{canonicalformmatrix}
  \begin{array}{r}
   \left. u_j = 1
   \begin{array}{l}
    \\
    \\
    \\
    \\
   \end{array}
    \right\{ \\
    \vdots \\
   \left. u_j = d_j
  \begin{array}{l}
    \\
    \\
    \\
    \\
    \end{array}
    \right\{ \\
   \left. u_k = 1
    \begin{array}{l}
    \\
    \\
    \\
    \\
    \end{array}
    \right\{ \\
    \vdots \\
    \left. u_k = d_k
    \begin{array}{l}
    \\
    \\
    \\
    \\
    \end{array}
    \right\{
  \end{array}
 \left[
\begin{array}{llllllllll}
  \tilde{\rho}_{s_j,s_k}\tilde{z}_{s_j}^1 & & & & - \tilde{\rho}_{s_i,s_k} \tilde{z}_{s_i}^1 & & & & & \\
  & \tilde{\rho}_{s_j,s_k}\tilde{z}_{s_j}^1 & & & - \tilde{\rho}_{s_i,s_k} \tilde{z}_{s_i}^2 & & & & & \\
  & & \ddots &  & \vdots & & & & & \\
  & & & \tilde{\rho}_{s_j,s_k} \tilde{z}_{s_j}^1 & - \tilde{\rho}_{s_i,s_k} \tilde{z}_{s_i}^{d_i} & & & & &\\
  & & \ddots  & & & \ddots & & & \\
 \tilde{\rho}_{s_j,s_k} \tilde{z}_{s_j}^{d_j} & & & & & & - \tilde{\rho}_{s_i,s_k} \tilde{z}_{s_i}^1 & & &\\
 &  \tilde{\rho}_{s_j,s_k} \tilde{z}_{s_j}^{d_j} & & & & & - \tilde{\rho}_{s_i,s_k} \tilde{z}_{s_i}^2 & & & \\
 & & \ddots & & & & \vdots & & & \\
 && & \tilde{\rho}_{s_j,s_k} \tilde{z}_{s_j}^{d_j} & & & - \tilde{\rho}_{s_i,s_k} \tilde{z}_{s_i}^{d_i} & & & \\
   \tilde{\rho}_{s_j,s_k}\tilde{z}_{s_k}^1 & & & & & & &-\tilde{\rho}_{s_i,s_j} \tilde{z}_{s_i}^1 & & \\
  & \tilde{\rho}_{s_j,s_k}\tilde{z}_{s_k}^1 & & & & & &- \tilde{\rho}_{s_i,s_j} \tilde{z}_{s_i}^2 & & \\
  & & \ddots &  & & & & \vdots & &\\
  & & & \tilde{\rho}_{s_j,s_k} \tilde{z}_{s_k}^1 & & & & - \tilde{\rho}_{s_i,s_j} \tilde{z}_{s_i}^{d_i} & &\\
  & & \ddots & & & & & & \ddots & \\
 \tilde{\rho}_{s_j,s_k} \tilde{z}_{s_k}^{d_k} & & & & & & & & & - \tilde{\rho}_{s_i,s_j} \tilde{z}_{s_i}^1 \\
 &  \tilde{\rho}_{s_j,s_k} \tilde{z}_{s_k}^{d_k}  & & & & & & & & -\tilde{\rho}_{s_i,s_j} \tilde{z}_{s_i}^2 \\
 & & \ddots & & & & & & & \vdots \\
 && & \tilde{\rho}_{s_j,s_k} \tilde{z}_{s_k}^{d_k} & & & & & & - \tilde{\rho}_{s_i,s_j} \tilde{z}_{s_i}^{d_i} \\
 \end{array}
 \right]
\end{equation}
for each given triple \(i, j\) and \(k\).

In Section \ref{se:Notation_Background}, \(M_n\) is introduced as a set containing \(n\) arbitrary elements. Hereafter, we let \(M_n = \{1,2,\cdots,n\}\). Denote the group operation of \(G_{M_n}\) composition  by \(\circ\), and denote the identity permutation of \(M_n\) by \(e\), e.g. \(e(i)=i\), \(\forall i\in M_n\). The generators of \(G_{M_n}\) written in cyclic form are \((i\hspace{2mm} i+1)\), \(\forall i\in M_n\setminus\{n\}\).

\begin{prop}\label{permutation:necesuffi} Assume \(n\ge 3\) and \(A=\{s_i: 1\le i\le n\}\subseteq S\) with all sites in \(A\) being neighbors of each other under \(\partial\). Then, a necessary and sufficient condition for the permutation property of the joint distribution of \(X_{s_1}, ..., X_{s_n}\), i.e. \(\Pi_{s_1, ..., s_n}(x_{s_{1}}, ..., x_{s_{n}})\) is that (\ref{blockwiseCanonicalformsolutionform1}) is a solution of linear equations of blockwise Canonical form (\ref{blockwiseCanonicalform}).
\end{prop}

\proof For a given permutation \(g\in G_{M_n}\), similar to (\ref{closedform_neighbors}), the joint pmf of \(X_{s_{g(1)}}, ..., X_{s_{g(n)}}\) is
\begin{eqnarray}\label{closedform_neighborsgroupg}
 \Pi_{s_{g(1)}, ..., s_{g(n)}}(x_{s_{g(1)}}, ..., x_{s_{g(n)}}) &=& \sum_{i=2}^{n} \biggl[\biggl(\tilde{z}_{s_{g(i)}}\sum_{j=1}^{i-1} \prod_{k=1,k\not=j}^{i-1}\hat{\pi}_{s_{g(k)}}(x_{s_{g(k)}}) \tilde{\rho}_{s_{g(i)},s_{g(j)}}\tilde{z}_{s_{g(j)}}\biggl) \cdot \prod_{k=i+1}^n \pi_{s_{g(k)}}(x_{s_{g(k)}}) \biggl] \nonumber \\
 &&  + \prod_{i=1}^n \pi_{s_{g(i)}}(x_{s_{g(i)}})
\end{eqnarray}
for each \(x_{s_i} \in \mathbf{X}_{s_i}\) (\(1\le i\le n\)).
A necessary and sufficient condition for the permutation property of the joint distribution of \(X_{s_1}, ..., X_{s_n}\) is that for any pair of permutations \(g, h\in G_{M_n}\),
\begin{equation}\label{conditionpermutation1}
\Pi_{s_{g(1)}, ..., s_{g(n)}}(x_{s_{g(1)}}, ..., x_{s_{g(n)}}) = \Pi_{s_{h(1)}, ..., s_{h(n)}}(x_{s_{h(1)}}, ..., x_{s_{h(n)}}).
\end{equation}
First, we prove necessity: For any \(g \in G_{M_n}\), there exists \(h\in G_{M_n}\) satisfying \(h(2) = g(3)\), \(h(3)=g(2)\) and \(h(i) = g(i)\), \(\forall i\in M_n\setminus\{2,3\}\) (We can also choose \(h\in G_{M_n}\) such that \(h(1)=g(3)\), \(h(2)=g(1)\), \(h(3)=g(2)\) and \(h(i) = g(i)\), \(\forall i\in M_n\setminus\{1,2,3\}\)). Then it follows by (\ref{closedform_neighborsgroupg}) and (\ref{conditionpermutation1}) that
\begin{equation}\label{necessarypermcond1}
  \tilde{\rho}_{s_{g(1)},s_{g(2)}}\tilde{z}_{s_{g(1)}}\tilde{z}_{s_{g(2)}} (\pi_{s_{g(3)}}(x_{s_{g(3)}}) - \hat{\pi}_{s_{g(3)}}(x_{s_{g(3)}})) = \tilde{\rho}_{s_{g(1)},s_{g(3)}}\tilde{z}_{s_{g(1)}}\tilde{z}_{s_{g(3)}} (\pi_{s_{g(2)}}(x_{s_{g(2)}}) - \hat{\pi}_{s_{g(2)}}(x_{s_{g(2)}})).
\end{equation}
Similarly, for the same \(g\), there exists \(h\in G_{M_n}\) satisfying \(h(1) = g(3)\), \(h(3)=g(1)\) and \(h(i) = g(i)\), \(\forall i\in M_n\setminus\{1,3\}\) (Another choice of \(h\in G_{M_n}\) is that \(h(1)=g(2)\), \(h(2)=g(3)\), \(h(3)=g(1)\) and \(h(i) = g(i)\), \(\forall i\in M_n\setminus\{1,2,3\}\)). Then we have that
\begin{equation}\label{necessarypermcond2}
  \tilde{\rho}_{s_{g(1)},s_{g(2)}}\tilde{z}_{s_{g(1)}}\tilde{z}_{s_{g(2)}} (\pi_{s_{g(3)}}(x_{s_{g(3)}}) - \hat{\pi}_{s_{g(3)}}(x_{s_{g(3)}})) = \tilde{\rho}_{s_{g(2)},s_{g(3)}}\tilde{z}_{s_{g(2)}}\tilde{z}_{s_{g(3)}} (\pi_{s_{g(1)}}(x_{s_{g(1)}}) - \hat{\pi}_{s_{g(1)}}(x_{s_{g(1)}})).
\end{equation}
Since \(g\) is arbitrary, it follows by equalities (\ref{necessarypermcond1}) and (\ref{necessarypermcond2}) that
\begin{equation}\label{necessarypermcond3}
  \tilde{\rho}_{s_{i},s_{j}}\tilde{z}_{s_{i}}\tilde{z}_{s_{j}} (\pi_{s_{k}}(x_{s_{k}}) - \hat{\pi}_{s_{k}}(x_{s_{k}}))
  = \tilde{\rho}_{s_{i},s_{k}}\tilde{z}_{s_{i}}\tilde{z}_{s_{k}} (\pi_{s_{j}}(x_{s_{j}}) - \hat{\pi}_{s_{j}}(x_{s_{j}}))
  = \tilde{\rho}_{s_{j},s_{k}}\tilde{z}_{s_{j}}\tilde{z}_{s_{k}} (\pi_{s_{i}}(x_{s_{i}}) - \hat{\pi}_{s_{i}}(x_{s_{i}}))
\end{equation}
for all \(x_{s_{i}}\in \mathbf{X}_{s_i}\), \(x_{s_{j}}\in \mathbf{X}_{s_j}\), \(x_{s_{k}}\in \mathbf{X}_{s_k}\) and distinct triple \(1\le i,j,k \le n\). (\ref{necessarypermcond3}) is equivalent to (\ref{blockwiseCanonicalformsolutionform1}) being a solution of linear equations of blockwise Canonical form (\ref{blockwiseCanonicalform}).

Now we prove sufficiency: If (\ref{blockwiseCanonicalformsolutionform1}) is a solution of linear equations of blockwise Canonical form (\ref{blockwiseCanonicalform}), then (\ref{necessarypermcond3}) follows. Next, we prove that for each \(g\in G_{M_n}\), \(\forall i\in M_n\setminus\{n\}\), we have
\begin{equation}\label{sufficonditionpermutation}
\Pi_{s_{g(1)}, ..., s_{g(n)}}(x_{s_{g(1)}}, ..., x_{s_{g(n)}}) = \Pi_{s_{g\circ (i\hspace{1mm} i+1)(1)}, ..., s_{g\circ (i\hspace{1mm} i+1)(n)}}(x_{s_{g\circ (i\hspace{1mm} i+1)(1)}}, ..., x_{s_{g\circ (i\hspace{1mm} i+1)(n)}}).
\end{equation}
By using (\ref{closedform_neighborsgroupg}), we can verify (\ref{sufficonditionpermutation}) holds for \(i=1\) without any other condition. For fixed \(2\le i\le n-1\),
(\ref{sufficonditionpermutation}) is equivalent to
\begin{eqnarray}\label{sufficonditionpermutation2}
 && \Pi_{s_{g(1)}, ..., s_{g(i-1)}, s_{g(i)}, s_{g(i+1)},..., s_{g(n)}}(x_{s_{g(1)}}, ..., x_{s_{g(i-1)}}, x_{s_{g(i)}}, x_{s_{g(i+1)}},..., x_{s_{g(n)}}) \nonumber \\
&=& \Pi_{s_{g(1)}, ..., s_{g(i-1)}, s_{g(i+1)}, s_{g(i)},..., s_{g(n)}}(x_{s_{g(1)}}, ..., x_{s_{g(i-1)}}, x_{s_{g(i+1)}}, x_{s_{g(i)}},..., x_{s_{g(n)}}),
\end{eqnarray}
which though (\ref{closedform_neighborsgroupg}) is guaranteed by
\begin{equation}\label{sufficonditionpermutation3}
  \tilde{\rho}_{s_{g(j)},s_{g(i)}}\tilde{z}_{s_{g(j)}}\tilde{z}_{s_{g(i)}} (\pi_{s_{g(i+1)}}(x_{s_{g(i+1)}}) - \hat{\pi}_{s_{g(i+1)}}(x_{s_{g(i+1)}}))
  = \tilde{\rho}_{s_{g(j)},s_{g(i+1)}}\tilde{z}_{s_{g(j)}}\tilde{z}_{s_{g(i+1)}} (\pi_{s_{g(i)}}(x_{s_{g(i)}}) - \hat{\pi}_{s_{g(i)}}(x_{s_{g(i)}})),
\end{equation}
\(\forall 1\le j\le i-1\). For each \(1\le j\le i-1\), \(2\le i\le n-1\) and \(g\in G_{M_n}\), (\ref{sufficonditionpermutation3}) follows by (\ref{necessarypermcond3}). Therefore, (\ref{sufficonditionpermutation}) is proved. Let \(g=e\), then (\ref{sufficonditionpermutation}) implies that for each generator \((i\hspace{2mm} i+1)\) (\(1\le i\le n-1\)),
\begin{equation}\label{sufficonditionpermutationgenerator}
\Pi_{s_{1}, ..., s_{n}}(x_{s_{1}}, ..., x_{s_{n}}) = \Pi_{s_{(i\hspace{1mm} i+1)(1)}, ..., s_{(i\hspace{1mm} i+1)(n)}}(x_{s_{(i\hspace{1mm} i+1)(1)}}, ..., x_{s_{(i\hspace{1mm} i+1)(n)}}).
\end{equation}
(\ref{sufficonditionpermutation}) and (\ref{sufficonditionpermutationgenerator}) together lead to the permutation property of the joint distribution of \(X_{s_1}, ..., X_{s_n}\). \qed

Before we discuss the solution space of the linear equations of blockwise Canonical form (\ref{blockwiseCanonicalform}), we point out its trivial solution \(\mathbf{0}\) where \(\mathbf{0}\) is the \(d\)-mensional zero vector. This \(\mathbf{0}\) implies a simple albeit special condition on the permutation property of \(X_{s_1}, ..., X_{s_n}\).

\begin{cor}\label{cor:permutation:necesuffi} Assume \(n\ge 3\) and \(A=\{s_i: 1\le i\le n\}\subseteq S\) with all sites in \(A\) being neighbors of each other under \(\partial\). Then, a sufficient condition for the permutation property of the joint distribution of \(X_{s_1}, ..., X_{s_n}\) is that \(\hat{\pi}_{s_i}(\cdot) \equiv \pi_{s_i}(\cdot)\) for \(1\le i\le n\).
\end{cor}
\proof For each \(1\le i\le n\) and \(1\le u_i\le d_i\), \(\hat{\pi}_{s_i}(\cdot) \equiv \pi_{s_i}(\cdot)\) implies that \(\hat{y}_{s_i}^{u_i}=\hat{\pi}_{s_i}(x_{s_i}^{u_i}) - \pi_{s_i}(x_{s_i}^{u_i})=0\), i.e. (\ref{blockwiseCanonicalformsolutionform1}) gives the \(\mathbf{0}\) solution of (\ref{blockwiseCanonicalform}). The conclusion of this corollary follows by the sufficiency of Proposition \ref{permutation:necesuffi}. \qed

We make a very loose assumption that no distribution in \(\{\pi_s\}_{s\in S}\), \(\{\tilde{\pi}_s\}_{s\in S}\) and \(\{\hat{\pi}_s\}_{s\in S}\) is degenerate, i.e., being the distribution of a constant random variable. Now we discuss the structure of the solution space of the linear equations of (\ref{blockwiseCanonicalform}). The solution space of (\ref{blockwiseCanonicalform}) is the \textbf{null} space of the canonical form matrix in (\ref{canonicalformmatrix}). This space depends on the value of \(n\) and the correlation structure among \(X_{s_1}, ..., X_{s_n}\), i.e. \(\rho_{s_i,s_j}\) (\(1\le i<j\le n\)). We deal with the case \(n=3\) and the case \(n\ge 4\) separately.

When \(n=3\), \(i,j\) and \(k\) essentially have only one combination, e.g., \(i=1, j=2\) and \(k=3\). Based on the assumption above, there exists at least one  \(1\le \check{u}_1 \le d_1\) such that \(\tilde{z}_{s_1}^{\check{u}_1}\not=0\). It follows that (\ref{blockwiseCanonicalform}) has solutions:
\begin{eqnarray}\label{n=3,generalsolution}
  y_{s_1}^{u_1} & = & \frac{\tilde{z}_{s_1}^{u_1} y_{s_1}^{\check{u}_1}}{\tilde{z}_{s_1}^{\check{u}_1}}  \hspace{2.8cm} \forall 1\le u_1 \le d_1  \nonumber \\
  y_{s_2}^{u_2} & = & \frac{\tilde{\rho}_{s_2,s_3} \tilde{z}_{s_2}^{u_2} y_{s_1}^{\check{u}_1}}{\tilde{\rho}_{s_1,s_3} \tilde{z}_{s_1}^{\check{u}_1}}  \hspace{2cm} \forall 1\le u_2 \le d_2  \\
  y_{s_3}^{u_3} & = & \frac{\tilde{\rho}_{s_2,s_3}\tilde{z}_{s_3}^{u_3} y_{s_1}^{\check{u}_1}}{\tilde{\rho}_{s_1,s_2}\tilde{z}_{s_1}^{\check{u}_1}}  \hspace{2cm} \forall 1\le u_3 \le d_3 \nonumber
\end{eqnarray}
where \(y_{s_1}^{\check{u}_1}\) is a free variable. It means that the solution space of (\ref{blockwiseCanonicalform}) is one-dimensional.

When \(n\ge 4\), it becomes more complicated. First, we define the correlation multiplication equalities.
{\defn\label{correlationmultiplicationeqncond} \textbf{Correlation multiplication equalities} are
\begin{equation}\label{correlationmultiplicationequalities}
    \rho_{s_i,s_j}\rho_{s_k,s_l} = \rho_{s_i,s_l}\rho_{s_j,s_k} = \rho_{s_i,s_k}\rho_{s_j,s_l}
\end{equation}
for each quadruple \(1\le i,j,k,l \le n\) where \(i,j,k\) and \(l\) are distinct integers.}

In Example \ref{exampleCorrMultiplicationEqn}, we will construct a random field whose nearby correlations indeed satisfy the correlation multiplication equalities (\ref{correlationmultiplicationequalities}).

For each quadruple \(1\le i,j,k,l \le n\) with \(i,j,k,l\) being distinct, we may consider the triples (\(i,j,k\)), (\(i,j,l\)), (\(i,k,l\)) and (\(j,k,l\)) and apply the \(n=3\) case to each of them. Each of the four groups of equations of blockwise canonical form has solutions formally similar to (\ref{n=3,generalsolution}), that is to say, we can solve
\(\{y_{s_i}^{u_i}, y_{s_j}^{u_j}, y_{s_k}^{u_k}\}_{(u_i, u_j, u_k)=(1,1,1)}^{(d_i, d_j,d_k)}\), \(\{y_{s_i}^{u_i}, y_{s_j}^{u_j}, y_{s_l}^{u_l}\}_{(u_i, u_j, u_l)=(1,1,1)}^{(d_i, d_j,d_l)}\), \(\{y_{s_i}^{u_i}, y_{s_k}^{u_k}, y_{s_l}^{u_l}\}_{(u_i, u_k, u_l)=(1,1,1)}^{(d_i, d_k,d_l)}\) and \newline \(\{y_{s_j}^{u_j}, y_{s_k}^{u_k}, y_{s_l}^{u_l}\}_{(u_j, u_k, u_l)=(1,1,1)}^{(d_j, d_k,d_l)}\) separately. Then same variables, e.g. \(\{y_{s_i}^{u_i}\}_{u_i=1}^{d_i}\) have solutions to three group of equations with unknown variables
\(\{y_{s_i}^{u_i}, y_{s_j}^{u_j}, y_{s_k}^{u_k}\}_{(u_i, u_j, u_k)=(1,1,1)}^{(d_i, d_j,d_k)}\), \(\{y_{s_i}^{u_i}, y_{s_j}^{u_j}, y_{s_l}^{u_l}\}_{(u_i, u_j, u_l)=(1,1,1)}^{(d_i, d_j,d_l)}\), and \(\{y_{s_i}^{u_i}, y_{s_k}^{u_k}, y_{s_l}^{u_l}\}_{(u_i, u_k, u_l)=(1,1,1)}^{(d_i, d_k,d_l)}\) respectively. It follows that (\ref{correlationmultiplicationequalities}) is a sufficient and necessary condition for  \(\{y_{s_i}^{u_i}\}_{u_i=1}^{d_i}\),  \(\{y_{s_j}^{u_j}\}_{u_j=1}^{d_j}\),  \(\{y_{s_k}^{u_k}\}_{u_k=1}^{d_k}\) and  \(\{y_{s_l}^{u_l}\}_{u_l=1}^{d_l}\) have coincident solutions to equations of blockwise canonical form involved. Again, we assume that \(\tilde{z}_{s_1}^{\check{u}_1}\not=0\) for some \(1\le \check{u}_1 \le d_1\). Then we can verify that the \textbf{correlation multiplication equalities} are equivalent to that (\ref{blockwiseCanonicalform}) has solutions:
\begin{eqnarray}\label{n>=4,generalsolution}
 y_{s_1}^{u_1} & = & \frac{\tilde{z}_{s_1}^{u_1} y_{s_1}^{\check{u}_1}}{\tilde{z}_{s_1}^{\check{u}_1}}  \hspace{2.8cm} \forall 1\le u_1 \le d_1 \nonumber \\
 y_{s_i}^{u_i} & = & \frac{\tilde{\rho}_{s_i,s_n} \tilde{z}_{s_i}^{u_i} y_{s_1}^{\check{u}_1}}{\tilde{\rho}_{s_1,s_n} \tilde{z}_{s_1}^{\check{u}_1}}  \hspace{2cm}  \forall 1\le u_i \le d_i, 2\le i\le n-1  \\
 y_{s_n}^{u_n} & = & \frac{\tilde{\rho}_{s_2,s_n}\tilde{z}_{s_n}^{u_n} y_{s_1}^{\check{u}_1}}{\tilde{\rho}_{s_1,s_2}\tilde{z}_{s_1}^{\check{u}_1}}  \hspace{1.9cm} \forall 1\le u_n \le d_n \nonumber
\end{eqnarray}
where \(y_{s_1}^{\check{u}_1}\) is a free variable. It means that for the case \(n\ge 4\), when the conditions (\ref{correlationmultiplicationequalities}) are satisfied, the solution space of (\ref{blockwiseCanonicalform}) is also one-dimensional.

To summarize, when \(n=3\), (\ref{blockwiseCanonicalform}) has nonzero solutions expressed by (\ref{n=3,generalsolution}). When \(n\ge 4\), (\ref{blockwiseCanonicalform}) also has nonzero solutions formulated in (\ref{n>=4,generalsolution}), given the \textbf{correlation multiplication equalities} conditions (\ref{correlationmultiplicationequalities}) are satisfied. The conditions (\ref{correlationmultiplicationequalities}) themselves are not simple though. We need to further our investigation on the conditions that can guarantee the existence of (\ref{correlationmultiplicationequalities}).

\begin{lem}\label{lem:conditionBlockwiseCanonicalform} Assume \(n\ge 4\), \({\rho}_{s_1, s_2}\not=0\) and \({\rho}_{s_1, s_n}\not=0\). The \textbf{correlation multiplication equalities} conditions (\ref{correlationmultiplicationequalities}) are satisfied if and only if
\begin{eqnarray}\label{conditionBlockwiseCanonicalform}
 {\rho}_{s_2,s_k} & = & \frac{{\rho}_{s_2,s_n}{\rho}_{s_1,s_k}}{{\rho}_{s_1,s_n}}  \hspace{2.8cm} \forall 3\le k \le n-1 \nonumber \\
 {\rho}_{s_i,s_j} & = & \frac{{\rho}_{s_1,s_i}{\rho}_{s_2,s_n}{\rho}_{s_1,s_j}}{{\rho}_{s_1,s_2}{\rho}_{s_1,s_n}} \hspace{2.0cm} \forall 3\le i\not=j \le n-1 \\
 {\rho}_{s_k,s_n} & = & \frac{{\rho}_{s_2,s_n}{\rho}_{s_1,s_k}}{{\rho}_{s_1,s_2}}  \hspace{2.8cm} \forall 3\le k \le n-1, \nonumber
\end{eqnarray}
where \(\rho_{s_1,s_k}\) (\(2\le k\le n\)) and \(\rho_{s_2,s_n}\) are \(n\) free variables.
\end{lem}

\proof In (\ref{correlationmultiplicationequalities}), there are \(\binom{n}{2}\) many  \(\rho\)\,s and \(2\binom{n}{4}\) equations. One can easily verify that (\ref{conditionBlockwiseCanonicalform}) solves (\ref{correlationmultiplicationequalities}). On the contrary, for given \(n\) free variables \(\rho_{s_1,s_k}\) (\(2\le k\le n\)) and \(\rho_{s_2,s_n}\), we solve (\ref{correlationmultiplicationequalities}) and obtain its only solution (\ref{conditionBlockwiseCanonicalform}) as follows. For \(3\le k \le n-1\), apply (\ref{correlationmultiplicationequalities}) to \(1,2,k,n\), then \(\rho_{s_1,s_2}\rho_{s_k,s_n} = \rho_{s_1,s_n}\rho_{s_2,s_k} = \rho_{s_1,s_k}\rho_{s_2,s_n}\) follows, which results in the first and third equalities of (\ref{conditionBlockwiseCanonicalform}). For \(3\le i\not=j \le n-1\), we apply (\ref{correlationmultiplicationequalities}) to \(1,2,i,j\). As a result, \(\rho_{s_1,s_2}\rho_{s_i,s_j} = \rho_{s_1,s_i}\rho_{s_2,s_j}\) holds. Then we plug in the first equality of (\ref{conditionBlockwiseCanonicalform}) and the second equality of (\ref{conditionBlockwiseCanonicalform}) follows. \qed

Proposition \ref{permutation:necesuffm} follows by  Proposition \ref{permutation:necesuffi}, (\ref{n=3,generalsolution}) and (\ref{n>=4,generalsolution}), and Lemma \ref{lem:conditionBlockwiseCanonicalform}:
\begin{prop}\label{permutation:necesuffm} Assume that \(\{s_1, ..., s_n\}\) are \(n\ge 3\) sites of \((S,\partial)\) which are neighbors of each other. Assume for \(n\ge 3\), \({\rho}_{s_1, s_2}\not=0\), \({\rho}_{s_1, s_n}\not=0\), \(\tilde{z}_{s_1}^{\check{u}_1}\not=0\) for some \(1\le \check{u}_1 \le d_1\) and for \(n\ge 4\), (\ref{conditionBlockwiseCanonicalform}) holds. Then, a necessary and sufficient condition for the permutation property of the joint distribution of \(X_{s_1}, ..., X_{s_n}\) is the solutions of (\ref{blockwiseCanonicalform}), i.e. \(\hat{y}_{s_i}^{u_i}\) (\(1\le u_i\le d_i, 1\le i\le n\)) given by (\ref{blockwiseCanonicalformsolutionform1}), have closed form as follows:
\begin{eqnarray}\label{n>=4,generalsolutionhaty}
 \hat{y}_{s_1}^{u_1} & = & \hat{\pi}_{s_1}(x_{s_1}^{u_1}) - \pi_{s_1}(x_{s_1}^{u_1}) = \frac{\tilde{z}_{s_1}^{u_1} y_{s_1}^{\check{u}_1}}{\tilde{z}_{s_1}^{\check{u}_1}}  \hspace{2.8cm} \forall 1\le u_1 \le d_1 \nonumber \\
 \hat{y}_{s_i}^{u_i} & = &\hat{\pi}_{s_i}(x_{s_i}^{u_i}) - \pi_{s_i}(x_{s_i}^{u_i}) = \frac{\tilde{\rho}_{s_i,s_n} \tilde{z}_{s_i}^{u_i} y_{s_1}^{\check{u}_1}}{\tilde{\rho}_{s_1,s_n} \tilde{z}_{s_1}^{\check{u}_1}}  \hspace{2.1cm}  \forall 1\le u_i \le d_i, 2\le i\le n-1  \\
 \hat{y}_{s_n}^{u_n} & = & \hat{\pi}_{s_n}(x_{s_n}^{u_n}) - \pi_{s_n}(x_{s_n}^{u_n}) =  \frac{\tilde{\rho}_{s_2,s_n}\tilde{z}_{s_n}^{u_n} y_{s_1}^{\check{u}_1}}{\tilde{\rho}_{s_1,s_2}\tilde{z}_{s_1}^{\check{u}_1}}  \hspace{1.8cm} \forall 1\le u_n \le d_n. \nonumber
\end{eqnarray}
\end{prop}

We give an example which illustrates the application of Proposition \ref{permutation:necesuffm}, Corollary \ref{cor:permutation:necesuffi}, Lemma \ref{lem:jointpmf} and Corollary \ref{lem:jointpmf} as follows.

{\examp \label{exampleCorrMultiplicationEqn} Let \(S=\{1,2,3,4\}\), \(s_i=i\)
for \(1\le i \le 4\) and \(\partial\) contain the following subsets
of \(S\): \(\partial(1)=\{2,3,4\}\), \(\partial(2)=\{1,3,4\}\),
\(\partial(3)=\{1,2,4\}\), \(\partial(4)=\{1,2,3\}\). Under this \(\partial\), all sites in \(S\) are neighbors of each other.

Let the common space of states for each site \(1\le i\le 4\) be \(\mathbf{X}=\{1,-1\}\); \(\tilde{\pi}_i(x_i)=\frac{1}{2}\)
\(\forall x_i\in \mathbf{X}\); \(\hat{\pi}_i(1)=\hat{p}\) and \(\hat{\pi}_i(-1)=1-\hat{p}\); \(\pi_i(1)=p\) and \(\pi_i(-1)=1-p\) where \(0<\hat{p}, p<1\). It follows that \(\mu_{\tilde{\pi}_i}=0\) and \(\sigma_{\tilde{\pi}_i}^2=1\). For convenience, let \(x_i^1 = 1\) and \(x_i^2 = -1\) (\(1\le i\le 4\)). Then by (\ref{defntildez}), \(\tilde{z}_i^{k}=\frac{1}{2}x_i^k\) for \(1\le i\le 4\), \(k=1,2\).
Let \(\beta_{12}=\beta_{13}=\beta_{23}=\beta_{14}=\beta_{24}=\beta_{34}=\beta\), i.e., all covariances between nearby sites are the same, then \(\rho_{12}=\rho_{13}=\rho_{23}=\rho_{14}=\rho_{24}=\rho_{34}=\rho\), and the \textbf{Correlation multiplication equalities} (\ref{correlationmultiplicationequalities}) hold. By (\ref{defntilderho}), \(\tilde{\rho}_{12}=\tilde{\rho}_{13}=\tilde{\rho}_{23}=\tilde{\rho}_{14}=\tilde{\rho}_{24}=\tilde{\rho}_{34}=\tilde{\rho}=\beta\) for any duple \(1\le i\not=j\le 4\). Consequently, all coefficients of linear equations (\ref{blockwiseCanonicalform}) are determined.

If \(\hat{p}=p\), the sufficient condition \(\hat{\pi}_{i}(\cdot) \equiv \pi_{i}(\cdot)\) \((1\le i\le 4)\) in Corollary \ref{cor:permutation:necesuffi} is satisfied, then the permutation property of \(\Pi_{1,2,3,4}(x_1,x_2,x_3,x_4)\) follows; otherwise, \(\hat{p}\not=p\) and it follows from Proposition \ref{permutation:necesuffm} that \(\Pi_{1,2,3,4}(x_1,x_2,x_3,x_4)\) is still permutable because the solution of (\ref{blockwiseCanonicalform}) given by (\ref{blockwiseCanonicalformsolutionform1}) is
\begin{eqnarray}\label{n>=4,generalsolutionhaty}
 \hat{y}_{i}^{1} & = & \hat{p} - p \hspace{2.8cm} \forall 1\le i\le 4 \nonumber \\
 \hat{y}_{i}^{2} & = & p - \hat{p} \hspace{2.8cm} \forall 1\le i\le 4, \nonumber
\end{eqnarray}
which exactly takes the closed form given by (\ref{n>=4,generalsolutionhaty}). The permutable joint pmf of \(X_1, ..., X_4\) follows directly from (\ref{closedform_neighbors}) as
\begin{eqnarray}\label{examplejointpmf}
&&\Pi_{1,2,3,4}(x_1,x_2,x_3,x_4) \\
&=& \prod_{i=1}^4 \pi_i(x_i) + \frac{1}{4}\beta [x_1 x_2\pi_3(x_3)\pi_4(x_4)+ x_1 x_3\hat{\pi}_2(x_2)\pi_4(x_4) + x_2 x_3\hat{\pi}_1(x_1)\pi_4(x_4) \nonumber\\
&& + x_1 x_4\hat{\pi}_2(x_2)\hat{\pi}_3(x_3) + x_2 x_4\hat{\pi}_1(x_1)\hat{\pi}_3(x_3) + x_3 x_4\hat{\pi}_1(x_1)\hat{\pi}_2(x_2)] \nonumber
\end{eqnarray}
for each \(x_i \in \mathbf{X}\) (\(1\le i\le 4\)).

The condition on \(\beta\) such that the right hand side of (\ref{examplejointpmf}) is within [0,1] can be given in explicit inequalities of \(\beta\). For illustration purpose, we only consider the special case of \(\hat{p}=p=\frac{1}{2}\). For this case, the necessary condition obtained from Corollary {lem:jointpmf} is \(\beta\in [-1,1]\). One can verify that the necessary and sufficient condition is actually \(\beta\in [-\frac{1}{6}, \frac{1}{2}]\).
}

By formula (\ref{closedform_neighbors}), for every pair \(1\le i\not=j\le n\), if \(\rho_{s_i,s_j}=0\), \(X_{s_i}\) and \(X_{s_j}\) are independent. More generally, we have the following corollary about the \textbf{independence of state on a site with uncorrelated states on its neighbors}:
{\cor\label{independenceofuncorrelatedneighbors} Suppose \(\hat{\pi}_{s_i}(\cdot) \equiv \pi_{s_i}(\cdot)\) (\(1\le i\le n\)). If for some \(1\le i\le n\), \(X_{s_i}\) is uncorrelated with rest \(X_{s_j}\) (\(1\le j\not=i\le n\)), then \(X_{s_i}\) is independent of \(X_{s_j}\) (\(1\le j\not=i\le n\)).
}
\proof By Corollary \ref{cor:permutation:necesuffi}, we only need to prove the statement if \(X_{s_n}\) is uncorrelated with each \(X_{s_j} (1\le j\le n-1)\), then \(X_{s_n}\) is independent of \(\{X_{s_j}\}_{j=1}^{n-1}\). This statement follows immediately by using formula (\ref{closedform_neighbors}) for the cases \(n\) and \(n-1\) respectively. \qed

{\rem\label{Independenceofexclusivesegments} Assume \(\hat{\pi}_{s_i}(\cdot) \equiv \pi_{s_i}(\cdot)\) (\(1\le i\le n\)). Even if \((X_{s_1}, X_{s_2})^T\) and \((X_{s_3}, X_{s_4})^T\) are uncorrelated, \((X_{s_1}, X_{s_2})^T\) and \((X_{s_3}, X_{s_4})^T\) are not necessarily independent. But Corollary \ref{independenceofuncorrelatedneighbors} does provide us a way to producing a \textbf{true} Markov random field with the remarkable properties of marginality and permutation. Recall the concept of connected space (\(S,\partial\)) in Section \ref{se:Notation_Background}. We extend the neighborhood system \(\partial\) to \(\partial'\) as follows: for each site \(s\in S\), let \(\partial'(s) = S\setminus\{s\}\) and \(\beta_{s,t}=0\), \(\forall t\in \partial'(s)\setminus\partial(s) = S\setminus (\partial(s)\bigcup\{s\})\). What's more, for each site \(s\in S\), because \(X_s\) is uncorrelated with \(\{X_t: t\in S\setminus (\partial(s)\bigcup\{s\})\}\), \(X_s\) is independent of \(\{X_t: t\in S\setminus (\partial(s)\bigcup\{s\})\}\) by Corollary \ref{independenceofuncorrelatedneighbors}. Now, let's return to the original neighborhood system \(\partial\). The previous conclusion infers that the state of each \(s\in S\) depends only on states of its neighbors \(t\in \partial(s)\), regardless of rest of the states. This exactly matches the definition of Markov random field in Definition \ref{defnRandomField}.
}

Before we state an alternative algorithm to generate our KNW Markov random field with those good properties of marginality and permutation, we update formula (\ref{CondProb}) with the sufficient condition for permutation, i.e., \(\hat{\pi}_{s}(\cdot) \equiv \pi_{s}(\cdot)\) (\(s\in S\)). For any ordered sequence \(\{s_i\}_{i=1}^N\) of sites in \(S\), it follows by formula (\ref{CondProb}) and Corollary \ref{independenceofuncorrelatedneighbors} that
\begin{eqnarray}\label{newformulaforPiAsi}
    &&\Pi(X_{s_i} = x_{s_i}|X_{s_{i-1}}=x_{s_{i-1}}, ...,X_{s_{1}}=x_{s_1}) \nonumber \\
 &=& \Pi(X_{s_i} = x_{s_i}|X_{A_{s_i}}=x_{A_{s_i}})  \\
 &=& \pi_{s_i}(x_{s_i}) + \frac{\tilde{\pi}_{s_i}(x_{s_i})(x_{s_i}-\mu_{\tilde{\pi}_{s_i}})}{\sigma^2_{\tilde{\pi}_{s_i}}}\sum_{t_i\in A_{s_i}}\biggl(\prod_{u_i\in A_{s_i}\setminus \{t_i\}}\pi_{u_i}(x_{u_i})\biggl) \cdot \frac{\beta_{s_i,t_i}\tilde{\pi}_{t_i}(x_{t_i})(x_{t_i}-\mu_{\tilde{\pi}_{t_i}})}{\sigma^2_{\tilde{\pi}_{t_i}}\Pi(X_{A_{s_i}}=x_{A_{s_i}})},  \nonumber
\end{eqnarray}
where \(A_{s_i}= \partial(s_i)\bigcap\{s_1,...,s_{i-1}\}\).

Now we propose the alternative algorithm to generate the KNW Markov random field with the desired marginality and permutation properties on \(S\):

Do for \(i=1,\dots, N\):
\begin{enumerate}
\item Choose arbitrary site \(s\in S\setminus\{s_1,..., s_{i-1}\}\) and denote it by \(s_i\). Then, let \(A_{s_i}= \partial(s_i)\bigcap\{s_1,...,s_{i-1}\}\), and compute \( \Pi(X_{A_{s_i}}=x_{A_{s_i}})\) by treating \(A_{s_i}\) as \(S\) and recursively repeating all the steps listed here.

\item Based on \(\Pi(X_{A_{s_i}}=x_{A_{s_i}})\), we compute
\( \Pi(\left. X_{s_i}=x_{s_i}^{u_i}\right| X_{A_{s_i}}=x_{A_{s_i}}) \) for \( 1\le u_i\le d_i \), using (\ref{newformulaforPiAsi}).

\item Generate a \([0,1]\)-uniform random variable \(U\). If
\[
\sum_{u_i=1}^{j-1} \Pi(\left. X_{s_i}=x_{s_i}^{u_i}\right| X_{A_{s_i}}=x_{A_{s_i}}) \leq U < \sum_{u_i=1}^{j} \Pi(\left. X_{s_i}=x_{s_i}^{u_i}\right| X_{A_{s_i}}=x_{A_{s_i}})
\]
 for some \( 1\le j\le d_i \), then set \(X_{s_i}= x_{s_i}^{j}\).
 We use \(x_{s_i}\) to indicate the simulated value \(x_{s_i}^{j}\) from \(\mathbf{X}_{s_i}\).
\end{enumerate}

{\rem\label{newrandomfield} Note that in the generation of random field \(X_S\) on (\(S,\partial\)), we made the extra assumption that for each site \(s\in S\), \(\beta_{s,t}=0\), \(\forall t\in S\setminus (\partial(s)\bigcup\{s\})\). The gain is the simplicity of the new algorithm compared to the algorithm in Section \ref{se:Algor_Sim_Ran_Field} and the desirable properties of marginality and permutation. This extra condition is actually necessary for generating a Markov random field, since state of each site on a Markov random field has to be independent of states outside of its neighborhood so as to be uncorrelated with them. Note also that \(S=\{s_i\}_{i=1}^N\) is not prescribed, but dynamically chosen. Here we took the great advantage of permutation which guarantees that our simulated Markov random field does not depend on how we order sites in \(S\). This alternative algorithm behaves in Gibbs sampler manner because every time it simulates state on a new site, it uses as much information as possible from simulated states in its neighborhood. But it is not Gibbs sampler algorithm, because this algorithm also simulates a Markov random field within one-pass.

Note also that in the above algorithm, \(\Pi(X_{A_{s_i}}=x_{A_{s_i}})\) is computed through treating \(A_{s_i}\) as \(S\) and recursively repeating all the steps in the algorithm, not through the multiplication rule employed in the algorithm of Section \ref{se:Algor_Sim_Ran_Field}. The reason is that the sites in \(A_{s_i}\) may not be neighbors of each other, if we use the multiplication rule, the Markov property to be established will be broken. Since there is no such concern as the Makrov property for the algorithm in Section \ref{se:Algor_Sim_Ran_Field}, multiplication rule can be applied there. Genearally speaking, \(S\) is big and neighborhoods are small, the method of computing \(\Pi(X_{A_{s_i}}=x_{A_{s_i}})\) in this alternative algorithm is much more efficient than that using multiplication rule.
}

\section{The Appendix}\label{se:Proofs_Lemma_Proposition}

\noindent {\em Proof of Lemma \ref{ConnectedNS}.} First, we prove that we can enumerate the sites within \(H\) in a sequence \(\{s_i\}_{i=1}^n\) such that \(s_i\in \partial(H^C\bigcup\{s_1,...,s_{i-1}\})\) for \(1\le i\le n\). Because \(H\) is nonempty, \(H^C\) is a proper subset of \(S\).  If \(H=S\),  \(\partial(H^C)=\partial(\emptyset)=S\), then we pick up arbitrary \(s_1\in \partial(H^C) = S\). Otherwise, \(H\) is a proper subset of \(S\), so \(H^C \not= \emptyset\). Because \(S\) is connected, \(\partial(H^C) \not= \emptyset\), we choose arbitrary \(s_1\in \partial(H^C)\). Assume that for some \(1\le j<n\) we have \(\{s_i\}_{i=1}^j\) such that \(s_i \in \partial(H^C\bigcup\{s_1,...,s_{i-1}\})\) for \(i=1, ..., j\). Then, we have \(s_i \notin H^C\bigcup \{s_1,...,s_{i-1}\}\) for \(i=1, ..., j\), which implies that \(s_i\) \((i=1, ..., j)\) are \(j\) different sites within \(H\). Because \(1\le j<n\) and \(H\) consists of \(n\) sites, \(H^C\bigcup\{s_1,...,s_{j}\}\) is a nonempty proper subset of \(S\). Since \(S\) is connected, \(\partial(H^C\bigcup\{s_1, ..., s_{j}\})\not=\emptyset\), we can find  some arbitrary \(s_{j+1}\in \partial(H^C\bigcup\{s_1, ..., s_{j}\})\). By induction, the statement is proved.

By the above construction, for each \(1\le i\le n\), \(s_i\in \partial(H^C\bigcup\{s_1,...,s_{i-1}\})\), it follows that \newline \(\partial (s_i)\bigcap [H^C\bigcup\{s_1,..., s_{i-1}\}]\not=\emptyset\) except when \(i=1\) and \(H=S\). If \(H=S\), then \(\partial (s_i)\bigcap [H^C\bigcup\{s_1,..., s_{i-1}\}]=\emptyset\), we let \(m_1=1\) and \(B_{s_1}^1=\emptyset\) such that \(\displaystyle \partial (s_i)\bigcap [H^C\bigcup\{s_1,..., s_{i-1}\}]=\bigcup_{j=1}^{m_i}B_{s_i}^j\). We also let \(A_{s_1}=B_{s_1}^1\). To prove the general case, it is sufficient for our purpose to prove for any nonempty \(A\subset S\), there exist \(m\ge 1\) and exclusive subsets \(\{B^j\}_{j=1}^{m}\) such that \(\displaystyle A=\bigcup_{j=1}^{m}B^j\) and each \(B^j\) (\(1\le j\le m\)) is connected.

For a given nonempty \(A\subset S\), we now present one way of choosing a nonempty connected component \(B\subseteq A\). We start \(B\) with the empty set (\(B=\emptyset\)). We select an arbitrary site \(t_1\in A\), and append \(t_1\) to \(B\) (\(B=\{t_1\}\)). If \(\partial(B)\bigcap A=\emptyset\), that is to say, \(B\) has no neighbors within \(A\), then we stop and \(B\) is finalized (\(B=\{t_1\}\)). If \(\partial(B)\bigcap A \not=\emptyset\), then we pick up a site \(t_2\in \partial(B)\bigcap A\), a neighbor of \(B\) within \(A\), and append \(t_2\) to \(B\) (\(B=\{t_1, t_2\}\)). If \(\partial(B)\bigcap A=\emptyset\), we stop and \(B\) is finalized (\(B=\{t_1, t_2\}\)). We repeat this procedure and add more and more neighbors of \(B\) within \(A\) to \(B\). Because \(A\) is finite, the procedure must be stopped for some \(1\le k \le |A|\) where \(|A|\) is the size of \(A\), then \(B\) is finalized (\(B=\{t_1, t_2, ..., t_k\}\)). By the construction of \(B\), \(t_j\in \partial(\{t_1, ..., t_{j-1}\})\) for \(2\le j\le k\), so \(B\) is nonempty and connected.

If \(B=A\), the claim in the second paragraph is proved. Otherwise, let \(B^1=B\). Then \(A\setminus B^1\not=\emptyset\), \(B^1\) and \(A\setminus B^1\) are exclusive. Then we can replace \(A\) in the third paragraph with \(A\setminus B^1\) and repeat the procedure there. Since the orginal \(A\) only has finite sites, there exist some \(m\ge 1\) for \(A\), such that we can repeat the procedure in the third paragraph \(m\) times to find out those \(m\) connected components of \(A\) denoted by \(\{B^j\}_{j=1}^m\) such that \(\displaystyle A=\bigcup_{j=1}^{m}B^j\) and \(\{B^j\}_{j=1}^m\) are mutually exclusive. For the given nonempty \(A\subset S\), once its nonempty exclusive connected components \(\{B^j\}_{j=1}^m\) are determined, we can always choose among them a component which contains the largest number of sites.

Henceforth, for each nonempty \(\partial (s_i)\bigcap [H^C\bigcup\{s_1,..., s_{i-1}\}]\), we first its exclusive connected components \(\{B_{s_i}^j\}_{j=1}^{m_i}\) such that \(\displaystyle \partial(s_i)\bigcap [H^C\bigcup\{s_1,..., s_{i-1}\}]=\bigcup_{j=1}^{m_i} B_{s_i}^j\) and then choose among \(\{B_{s_i}^j\}_{j=1}^{m_i}\) one of the largest (meaning the number of sites) and denote it by \(A_{s_i}\).
\endproof

\noindent {\em Proof of Proposition \ref{mainprop:alg}.} It is clear by (\ref{CondProb}) that for $i=1$, $X_{s_1}$ has probability distribution $\pi_{s_1}(\cdot)$, since $A_{s_1}=\emptyset$.

We next to prove that for $2\le i\le N$, $X_{s_i}$ has probability distribution $\pi_{s_i}(\cdot)$, and for any $t_i\in A_{s_i}$, $cov(X_{s_i}, X_{t_i})=\beta_{s_i, t_i}$. To ease notation, we suppress the subscript $i$. For \(x_s\in \mathbf{X}_s\), by (\ref{CondProb}), one has that
\begin{eqnarray}
\Pi(X_{s}=x_{s})
&=& \sum_{x_{A_s}\in \mathbf{X}_{A_s}} \Pi(\left. X_{s}=x_{s}\right| X_{A_s}=x_{A_s}) \Pi(X_{A_s}=x_{A_s}) \nonumber \\
&=& \pi_{s}(x_{s})+ \frac{\displaystyle{\tilde{\pi}_{s}(x_{s})(x_{s}-\mu_{\tilde{\pi}_{s}})}}
{\displaystyle{\sigma^2_{\tilde{\pi}_{s}}}}\sum_{t\in A_{s}} \sum_{x_{A_s}\in \mathbf{X}_{A_s}}
\bigg(\prod_{u\in A_{s}\setminus \{t\}}\hat{\pi}_{u}(x_{u})\cdot\frac{\displaystyle{
\beta_{s,t}\tilde{\pi}_{t}(x_{t})(x_{t}-\mu_{\tilde{\pi}_{t}})}}
{\displaystyle{\sigma^2_{\tilde{\pi}_{t}}}}\biggl) \nonumber \\
&=& \pi_{s}(x_{s})+ \frac{\displaystyle{\tilde{\pi}_{s}(x_{s})(x_{s}-\mu_{\tilde{\pi}_{s}})}}
{\displaystyle{\sigma^2_{\tilde{\pi}_{s}}}}\sum_{t\in A_{s}} \beta_{s,t} \sum_{x_{{A_s}\setminus\{t\}}\in \mathbf{X}_{{A_s}\setminus\{t\}}}\bigg(\prod_{u\in A_{s}\setminus \{t\}}\hat{\pi}_{u}(x_{u})\cdot \sum_{x_{t}\in \mathbf{X}_{t}} \frac{\displaystyle{\tilde{\pi}_{t}(x_{t})(x_{t}-\mu_{\tilde{\pi}_{t}})}}{\displaystyle{\sigma^2_{\tilde{\pi}_{t}}}}\biggl)\nonumber \\
&=& \pi_{s}(x_{s}), \label{marignal}
\end{eqnarray}
where the third equality follows from the equality \(\sum_{x_{A_s}\in \mathbf{X}_{A_s}}=\sum_{x_{{A_s}\setminus\{t\}}\in \mathbf{X}_{{A_s}\setminus\{t\}}}\sum_{x_{t}\in \mathbf{X}_{t}}\) and interchanging the order of summations, and the fourth equality holds, because for fixed $t\in {A_s}$,
$$
\sum_{x_{t}\in \mathbf{X}_{t}} \tilde{\pi}_{t}(x_{t})(x_{t}-\mu_{\tilde{\pi}_{t}})=\mu_{\tilde{\pi}_{t}} - \mu_{\tilde{\pi}_{t}} = 0.
$$

Now fix $t\in A_{s}$, we prove $cov(X_{s}, X_{t})=\beta_{s, t}$. We compute the joint probability mass function of $X_s$ and $X_t$. For  \(x_s \in \mathbf{X}_s, x_t\in \mathbf{X}_t\), we have that by (\ref{CondProb}) again
\begin{eqnarray}\label{Jointprob0}
&{}&\Pi(X_{s}=x_{s},  X_{t}=x_{t}) \nonumber \\
&=& \sum_{x_{{A_s}\setminus\{t\}}\in\mathbf{X}_{{A_s}\setminus\{t\}}} \Pi( X_{s}=x_{s}| X_{t}=x_{t}, X_{{A_s}\setminus\{t\}}=x_{{A_s}\setminus\{t\}})
\times \Pi( X_{t}=x_{t}, X_{{A_s}\setminus\{t\}}=x_{{A_s}\setminus\{t\}}) \nonumber \\
&=&\pi_{s}(x_{s})\sum_{x_{{A_s}\setminus\{t\}}\in\mathbf{X}_{{A_s}\setminus\{t\}}} \Pi( X_{t}=x_{t}, X_{{A_s}\setminus\{t\}}=x_{{A_s}\setminus\{t\}}) \nonumber \\
&{}& +  \frac{\displaystyle{\tilde{\pi}_{s}(x_{s})(x_{s}-\mu_{\tilde{\pi}_{s}})}}
{\displaystyle{\sigma^2_{\tilde{\pi}_{s}}}}\sum_{x_{{A_s}\setminus\{t\}}\in\mathbf{X}_{{A_s}\setminus\{t\}}}\sum_{u\in A_{s}} \biggl(\prod_{v\in A_{s}\setminus \{u\}}\hat{\pi}_{v}(x_{v}) \cdot\frac{\displaystyle{
\beta_{s,u}\tilde{\pi}_{u}(x_{u})(x_{u}-\mu_{\tilde{\pi}_{u}})}}
{\displaystyle{\sigma^2_{\tilde{\pi}_{u}}}}\biggl), \nonumber \\
\end{eqnarray}
since \( \Pi(X_{t}=x_{t}, X_{{A_s}\setminus\{t\}}=x_{{A_s}\setminus\{t\}}) = \Pi(X_{A_s}=x_{A_s}) \) for \(t\in A_s\). Therefore,
\begin{eqnarray}\label{Jointprob1}
&{}&\Pi(X_{s}=x_{s},  X_{t}=x_{t}) \nonumber \\
&=&\pi_{s}(x_{s})\pi_{t}(x_{t}) + \frac{\displaystyle{\tilde{\pi}_{s}(x_{s})(x_{s}-\mu_{\tilde{\pi}_{s}})}}
{\displaystyle{\sigma^2_{\tilde{\pi}_{s}}}}\sum_{x_{{A_s}\setminus\{t\}}\in\mathbf{X}_{{A_s}\setminus\{t\}}}\biggl[\prod_{v\in A_{s}\setminus \{t\}}\hat{\pi}_{v}(x_{v}) \cdot\frac{\displaystyle{
\beta_{s,t}\tilde{\pi}_{t}(x_{t})(x_{t}-\mu_{\tilde{\pi}_{t}})}}
{\displaystyle{\sigma^2_{\tilde{\pi}_{t}}}} \nonumber \\
&& + \sum_{u\in A_{s}\setminus\{t\}} \biggl(\prod_{v\in A_{s}\setminus \{u\}}\hat{\pi}_{v}(x_{v}) \cdot\frac{\displaystyle{
\beta_{s,u}\tilde{\pi}_{u}(x_{u})(x_{u}-\mu_{\tilde{\pi}_{u}})}}
{\displaystyle{\sigma^2_{\tilde{\pi}_{u}}}}\biggl)\biggl]   \nonumber \\
&=&\pi_{s}(x_{s})\pi_{t}(x_{t}) + \frac{\displaystyle{\tilde{\pi}_{s}(x_{s})(x_{s}-\mu_{\tilde{\pi}_{s}})\beta_{s,t}\tilde{\pi}_{t}(x_{t})(x_{t}-\mu_{\tilde{\pi}_{t}})}}
{\displaystyle{\sigma^2_{\tilde{\pi}_{s}}\displaystyle{\sigma^2_{\tilde{\pi}_{t}}}}}\cdot \prod_{v\in A_{s}\setminus \{t\}}\biggl(\sum_{x_{v}\in \mathbf{X}_{v}}\hat{\pi}_{v}(x_{v})\biggl) \nonumber \\
&& + \frac{\displaystyle{\tilde{\pi}_{s}(x_{s})(x_{s}-\mu_{\tilde{\pi}_{s}})}}
{\displaystyle{\sigma^2_{\tilde{\pi}_{s}}}}\sum_{u\in A_{s}\setminus\{t\}}\sum_{x_{{A_s}\setminus\{t\}}\in\mathbf{X}_{{A_s}\setminus\{t\}}} \biggl(\prod_{v\in A_{s}\setminus \{u\}}\hat{\pi}_{v}(x_{v}) \cdot\frac{\displaystyle{
\beta_{s,u}\tilde{\pi}_{u}(x_{u})(x_{u}-\mu_{\tilde{\pi}_{u}})}}
{\displaystyle{\sigma^2_{\tilde{\pi}_{u}}}}\biggl) \nonumber \\
&=&\pi_{s}(x_{s})\pi_{t}(x_{t}) + \frac{\displaystyle{\tilde{\pi}_{s}(x_{s})(x_{s}-\mu_{\tilde{\pi}_{s}})\beta_{s,t}\tilde{\pi}_{t}(x_{t})(x_{t}-\mu_{\tilde{\pi}_{t}})}}
{\displaystyle{\sigma^2_{\tilde{\pi}_{s}}\displaystyle{\sigma^2_{\tilde{\pi}_{t}}}}}
\end{eqnarray}
where the third equality follows from arguments similar to those in (\ref{marignal}), in particular:
\begin{eqnarray}
&&\sum_{u\in A_{s}\setminus\{t\}}\sum_{x_{{A_s}\setminus\{t\}}\in\mathbf{X}_{{A_s}\setminus\{t\}}} \biggl(\prod_{v\in A_{s}\setminus \{u\}}\hat{\pi}_{v}(x_{v}) \cdot\frac{\displaystyle{\beta_{s,u}\tilde{\pi}_{u}(x_{u})(x_{u}-\mu_{\tilde{\pi}_{u}})}}
{\displaystyle{\sigma^2_{\tilde{\pi}_{u}}}}\biggl) \nonumber \\
&=&\sum_{u\in A_{s}\setminus\{t\}}\biggl[\beta_{s,u} \sum_{x_{{A_s}\setminus\{t,u\}}\in\mathbf{X}_{{A_s}\setminus\{t,u\}}} \biggl(\hat{\pi}_{t}(x_{t}) \cdot \prod_{v\in A_{s}\setminus \{t,u\}}\hat{\pi}_{v}(x_{v})\biggl)\cdot\frac{\displaystyle{\sum_{x_{u}\in\mathbf{X}_{u}}}\tilde{\pi}_{u}(x_{u})(x_{u}-\mu_{\tilde{\pi}_{u}})}{\displaystyle{\sigma^2_{\tilde{\pi}_{u}}}}\biggl] \nonumber \\
&=&\hat{\pi}_{t}(x_{t})\sum_{u\in A_{s}\setminus\{t\}}\biggl[\beta_{s,u} \prod_{v\in A_{s}\setminus \{t,u\}} \biggl(\sum_{x_{v}\in\mathbf{X}_{v}}\hat{\pi}_{v}(x_{v})\biggl)\cdot\frac{\displaystyle{\sum_{x_{u}\in\mathbf{X}_{u}}}\tilde{\pi}_{u}(x_{u})(x_{u}-\mu_{\tilde{\pi}_{u}})}{\displaystyle{\sigma^2_{\tilde{\pi}_{u}}}}\biggl] \nonumber \\
&=&\hat{\pi}_{t}(x_{t})\sum_{u\in A_{s}\setminus\{t\}}\beta_{s,u}\cdot 0 \nonumber \\
&=& 0.
\end{eqnarray}

Therefore, by (\ref{Jointprob1}), we obtain,
\begin{eqnarray}\label{Jointprob2}
E[(X_{s}-\mu_{\tilde{\pi}_{s}})(X_{t}-\mu_{\tilde{\pi}_{t}})]
&=& \sum_{x_s\in \mathbf{X}_s}\sum_{x_t\in \mathbf{X}_t} (x_{s}-\mu_{\tilde{\pi}_{s}})(x_{t}-\mu_{\tilde{\pi}_{t}})\pi_{s}(x_{s})\pi_{t}(x_{t}) \nonumber \\
&& + \frac{\sum_{x_s\in \mathbf{X}_s}\sum_{x_t\in \mathbf{X}_t} \tilde{\pi}_{s}(x_{s})(x_{s}-\mu_{\tilde{\pi}_{s}})^2 \tilde{\pi}_{t}(x_{t})(x_{t}-\mu_{\tilde{\pi}_{t}})^2 \beta_{s,t}}
{\displaystyle{\sigma^2_{\tilde{\pi}_{s}}\displaystyle{\sigma^2_{\tilde{\pi}_{t}}}}}  \nonumber \\
&=&(E[X_{s}]-\mu_{\tilde{\pi}_{s}})(E[X_{t}]-\mu_{\tilde{\pi}_{t}}) + \beta_{s,t}.  \nonumber \\
\end{eqnarray}
It follows that
\begin{equation}\label{covarianceequality}
cov(X_{s},X_{t}) = cov(X_{s}-\mu_{\tilde{\pi}_{s}}, X_{t}-\mu_{\tilde{\pi}_{t}})=\beta_{s,t}
\end{equation} \qed
\endproof

{\rem Note that (\ref{Jointprob}) follows by (\ref{Jointprob1}) and (\ref{covarianceequality}).
}

\end{document}